\documentclass[11pt]{amsart}
\usepackage{amsxtra}
\usepackage{amssymb}
\theoremstyle{plain}

\newcommand{\cleqn}{\setcounter{equation}{0}}
\newcommand{\clth}{\setcounter{theorem}{0}}
\newcommand {\sectionnew}[1]{\section{#1}\cleqn\clth}
\newcommand{\nn}{\hfill\nonumber}
\newtheorem{theorem}{Theorem}[section]
\newtheorem{lemma}[theorem]{Lemma}
\newtheorem{definition-theorem}[theorem]{Definition-Theorem}
\newtheorem{proposition}[theorem]{Proposition}
\newtheorem{corollary}[theorem]{Corollary}
\newtheorem{definition}[theorem]{Definition}
\newtheorem{example}[theorem]{Example}
\newtheorem{remark}[theorem]{Remark}
\newtheorem{conjecture}[theorem]{Conjecture}
\newtheorem{notation}[theorem]{Notation}
\newcommand \bth[1] { \begin{theorem}\label{t#1} }
\newcommand \ble[1] { \begin{lemma}\label{l#1} }

\newcommand \bpr[1] { \begin{proposition}\label{p#1} }
\newcommand \bco[1] { \begin{corollary}\label{c#1} }
\newcommand \bde[1] { \begin{definition}\label{d#1}\rm }
\newcommand \bex[1] { \begin{example}\label{e#1}\rm }
\newcommand \bre[1] { \begin{remark}\label{r#1}\rm }
\newcommand \bcj[1] { \begin{conjecture}\label{j#1}\rm }

\newcommand \bnota[1] { \begin{notation}\label{n#1}\rm }
\renewcommand {\eth} { \end{theorem} }
\newcommand {\ele} { \end{lemma} }

\newcommand {\epr} { \end{proposition} }
\newcommand {\eco} { \end{corollary} }
\newcommand {\ede} { \end{definition} }
\newcommand {\eex} { \end{example} }
\newcommand {\ere} { \end{remark} }
\newcommand {\ecj} { \end{conjecture} }

\newcommand {\enota} { \end{notation} }
\newcommand \thref[1]{Theorem \ref{t#1}}
\newcommand \leref[1]{Lemma \ref{l#1}}
\newcommand \prref[1]{Proposition \ref{p#1}}

\newcommand \cjref[1]{Conjecture \ref{j#1}}
\newcommand \deref[1]{Definition \ref{d#1}}
\newcommand \exref[1]{Example \ref{e#1}}
\newcommand \reref[1]{Remark \ref{r#1}}
\newcommand \lb[1]{\label{#1}}

\def \Rset {{\mathbb R}}         
\def \Rp   {{\mathbb R}_{\geq 0}}
\def \Cset {{\mathbb C}}
\def \KK {{\mathbb K}}
\def \Zset {{\mathbb Z}}
\def \ZN   {{\mathbb Z}^{N}}
\def \ZNp  {{\mathbb Z}^N_+}
\def \Zp   {{\mathbb N}}
\def \Qset {{\mathbb Q}}

\def \Tset {{\mathbb T}}

\def \B  {{\mathcal{B}}}               
\def \AA {{\mathcal{A}}}
\def \II {{\mathcal{I}}}

\def \QQ {{\mathcal{Q}}}
\def \PP {{\mathcal{P}}}

\def \UU {{\mathcal{U}}}
\def \RR {{\mathcal{R}}}

\def \CC {{\mathcal{C}}}

\def \NN {{\mathcal{N}}}
\def \FF {{\mathcal{F}}}
\def \II {{\mathcal{I}}}
\def \HH {{\mathcal{H}}}
\def \TT {{\mathcal{T}}} 

\def \vv {{\nu}}
\def \ib {{\bf{i}}}
\def \qb {{\bf{q}}}
\def \db {{\bf{d}}}
\def \pb {{\bf{p}}}
\def \rb {{\bf{r}}}
\def \De {\Delta}   
\def \de {\delta}
\def \al {\alpha}
\def \be {\beta}
\def \vpi {\varpi}
\def \la {\lambda}

\def \om {\omega}
\def \Om {\Omega}
\def \ga {\gamma}
\def \de {\delta}
\def \Ga {\Gamma}

\def \sig {\sigma}

\def \sig{\sigma}

\def \mt  {\mapsto}

\def \hra {\hookrightarrow}

\def \ci  {\circ}

\def \rcor {\rangle}
\def \lcor {\langle}

\def \ol {\overline}
\def \wt {\widetilde}
\def \wh {\widehat}

\def \id { {\mathrm{id}} }

\def \sign { {\mathrm{sign}} }

\def \rank { {\mathrm{rank}} }

\def \g  {\mathfrak{g}}   

\def \sl {\mathfrak{sl}} 

\def \so {\mathfrak{so}}

\def \n  {\mathfrak{n}}

\def \sl {\mathfrak{sl}}

\DeclareMathOperator \Span { {\mathrm{Span}} }
\DeclareMathOperator \Aut { {\mathrm{Aut}} }

\DeclareMathOperator \Con { {\mathrm{Con}} }
\DeclareMathOperator \charr { {\mathrm{char}} }

\DeclareMathOperator \ad { {\mathrm{ad}} }

\DeclareMathOperator \Ker { {\mathrm{Ker}} }

\DeclareMathOperator \GKdim {{\mathrm{GK \, dim}}}

\DeclareMathOperator \Supp { {\mathrm{Supp}} }
\DeclareMathOperator \Fract { {\mathrm{Fract}} }

\renewcommand \max { {\mathrm{max}} }
\newcommand \Spec { {\mathrm{Spec}} }

\begin{document}
\title[Rigidity of quantum tori and the Andruskiewitsch--Dumas conjecture]
{Rigidity of quantum tori and the Andruskiewitsch--Dumas conjecture}
\author[Milen Yakimov]{Milen Yakimov}
\address{
Department of Mathematics \\
Louisiana State University \\
Baton Rouge, LA 70803
U.S.A.
}
\email{yakimov@math.lsu.edu}
\date{}
\begin{abstract} We prove the Andruskiewitsch--Dumas conjecture
that the automorphism group of the positive part of the quantized universal 
enveloping algebra $\UU_q(\g)$ of an arbitrary finite dimensional 
simple Lie algebra $\g$ is isomorphic to the semidirect product of the 
automorphism group of the Dynkin diagram of $\g$ and a 
torus of rank equal to the rank of $\g$. 

The key step in our proof is a rigidity theorem for quantum tori. 
It has a broad range of applications. It allows one to control the 
(full) automorphism groups of large classes of associative algebras, 
for instance quantum cluster algebras.
\end{abstract}
\maketitle
\sectionnew{Introduction}
\lb{intro}
Automorphism groups of algebras are often difficult to describe
and contain wild automorphisms. The latter fact was proved 
by Joseph \cite{J0} for $\Aut \UU(\sl_2)$, Alev \cite{A} for 
$\Aut \UU(\n)$ where $\n$ is the nilradical of a Borel subalgebra
of $\sl_3$, and Shestakov and Umirbaev \cite{SU} for the Nagata automorphism 
of a polynomial algebra in three variables.  

The Andruskiewitsch--Dumas conjecture \cite{AD} concerns the explicit
structure of the automorphism groups of the quantum nilpotent 
algebras $\UU_q^+(\g)$ for all simple Lie algebras $\g$. 
It asserts that, in contrast to the above cases, the algebras
$\UU_q^+(\g)$ are rigid in the sense that they have 
small automorphism groups. Despite many attempts to prove the conjecture, 
it remained wide open for $\g \neq \sl_3, \sl_4, \so_5$.
In this paper we prove the conjecture in full generality.

Let $\UU_q(\g)$ be the quantized universal enveloping algebra 
of a simple Lie algebra $\g$, defined over an arbitrary base field 
$\KK$ for a deformation parameter $q \in \KK^*$ which is not 
a root of unity. It has Cartan generators $E_\al$, $F_\al$, and 
$K_\al^{\pm 1}$, where $\al$ runs over the set $\Pi$ of simple roots of $\g$.
The algebra $\UU^+_q(\g)$ is the subalgebra of $\UU_q(\g)$ 
generated by $\{ E_\al \mid \al \in \Pi \}$. It is abstractly described as 
the $\KK$-algebra with those generators subject to  
the quantum Serre relations, see \eqref{S1}. The torus 
$\Tset^r = (\KK^*)^r$, 
where $r = |\Pi|$ is the rank of $\g$, 
acts on $\UU_q^+(\g)$ by algebra automorphisms by 
\[
t \cdot E_\al = t_\al E_\al, 
\; \; t = (t_{\al'})_{\al' \in \Pi}, 
\al \in \Pi. 
\]
The automorphism group of the Dynkin diagram $\Ga$ of $\g$ has a natural 
embedding into $\Aut(\UU^+_q(\g))$. To $\theta \in \Aut (\Ga)$ 
one associates the automorphism given by
\[
E_\al \mt E_{\theta(\al)},  
\; \; \al \in \Pi.
\]
Andruskiewitsch and Dumas \cite{AD} have conjectured that the above generate
the automorphism group $\Aut (\UU^+_q(\g))$.

\bcj{AD} (Andruskiewitsch--Dumas) For all simple Lie algebras $\g$ of rank $r>1$, 
base fields $\KK$, and deformation parameters $q \in \KK^*$ which are 
not roots of unity 
\[
\Aut (\UU_q^+(\g)) \cong \Tset^r \ltimes \Aut (\Ga).
\]
\ecj

Three cases of this conjecture were proved up to date: $\g = \sl_3$ 
by Alev--Dumas and Caldero \cite{AlD,Cal2}, 
$\g= \so_5$ by Launois \cite{Lau1} and Andruskiewitsch--Dumas \cite{AD},
and $\g = \sl_4$ by Lopes--Launois \cite{LL}. They
found important ways to study the automorphisms of $\UU_q^+(\g)$ from the 
structure of the spectra $\Spec \, \UU^+_q(\g)$. Unfortunately this relation 
could not be used to obtain sufficient restrictions on the automorphisms 
of $\UU^+_q(\g)$. As a result of this, the final steps of the proofs of 
the special cases relied on elaborate computations, which were specific to 
each case. Alev--Chamarie \cite{ACh}, Fleury \cite{F}, 
Launois--Lenagan \cite{LaL1,LaL2}, and Rigal \cite{R} studied 
the automorphism groups of quantum matrices, quantized universal enveloping 
algebras of Borel subalgebras, and quantum Weyl algebras.
In their works, arguments with induced actions on prime spectra and relations 
to derivations of quantum tori lead to enough information for automorphisms 
only when there were few height one primes, lots of units, or when the 
algebras had low Gelfand--Kirillov dimension. 

We give a proof of \cjref{AD} and exhibit a general 
classification method for automorphism groups of related algebras.
The key components of this new classification scheme are a relationship between 
$\Aut(\UU^+_q(\g))$ and the group of certain continuous bi-finite 
automorphisms of completed quantum tori, and a rigidity result for the latter. 
In order to state those, we need to introduce some more terminology and notation. 
Denote by $M_N(\KK^*)$ the set of $N \times N$ matrices with 
entries in $\KK^*$. An $N \times N$ matrix $\qb = (q_{kl})_{k,l=1}^N \in M_N(\KK^*)$
is called multiplicatively skew-symmetric if $q_{kl} q_{lk} = 1$
for $k \neq l$ and $q_{ll}= 1$. Such gives rise to the rank $N$ quantum torus
\begin{equation}
\label{q-tor}
\TT_\qb = \frac{\KK \lcor X_1^{\pm 1}, \ldots, X_N^{\pm 1} \rcor}{ (
X_k X_k^{-1} -1, X_k^{-1} X_k -1 , X_k X_l - q_{kl} X_l X_k, 1 \leq k < l \leq N )} \cdot
\end{equation}
Denote the multiplicative kernel of the matrix $\qb$
\begin{equation}
\label{Kerq}
\Ker (\qb) = \Big\{ (j_1, \ldots, j_N) \in \ZN
\mid \prod_{l=1}^N q_{kl}^{j_l} = 1, \; 
\forall 1 \leq k \leq N \Big\}. 
\end{equation}
Let $\Zset_+ :=\{1,2, \ldots\}$.
We say the quantum torus $\TT_\qb$ is {\em{saturated}} if 
\begin{equation}
\label{torsion-free}
f \in \ZN, n \in \Zset_+, \; \; 
n f \in \Ker(\qb) \Rightarrow f \in \Ker(\qb).
\end{equation}
For example, $\TT_\qb$ is saturated if the subgroup 
of $\KK^*$ generated by $q_{kl}$, 
$1 \leq k < l \leq N$ is torsion-free. 
The condition \eqref{torsion-free} has several other equivalent 
formulations, see \S \ref{3.1}. It is equivalent to the condition that for 
$u \in \TT_\qb$, $n \in \Zset_+$, $u^n \in Z(\TT_\qb)$ 
implies $u \in Z(\TT_\qb)$. Here and below for an 
algebra $B$, $Z(B)$ denotes its center. 
We call an $N$-tuple $\db = (d_1, \ldots, d_N) \in \Zset_+^N$
a degree vector and use it to define a $\Zset$-grading on 
$\TT_\qb$ by assigning $\deg X_l = d_l$. Consider the completion
\[
\wh{\TT}_{\qb, \db} := \{ u_m + u_{m+1} + \ldots \mid m \in \Zset,
u_j \in \TT_\qb, \deg u_j = j \}.
\] 
We will call a continuous automorphism $\phi$ of $\wh{\TT}_{\qb, \db}$ {\em{unipotent}} 
if 
\[
\phi(X_l) - X_l \in (\wh{\TT}_{\qb, \db})^{\geq d_l + 1}, \; \; \forall 1 \leq l \leq N,
\]
where $(\wh{\TT}_{\qb, \db})^{\geq m} := \{ u_m + u_{m+1} + \ldots \mid
u_j \in \TT_\qb, \deg u_j = j \}$ for $m \in \Zset$.
The $N$-tuple $(\phi(X_1), \ldots, \phi(X_N))$
consists of units of $\wh{\TT}_{\qb, \db}$ and 
uniquely determines the continuous automorphism $\phi$.
A unipotent automorphism $\phi$ of $\wh{\TT}_{\qb, \db}$ 
will be called {\em{bi-finite}} if 
\[
\phi(X_l) \; \mbox{and} \; \phi^{-1}(X_l) \in \TT_\qb, 
\; \; \forall 1 \leq l \leq N.
\]
We refer the reader to \S \ref{3.1} 
for properties of the above types of automorphisms.
To this end we note that in general bi-finite unipotent automorphisms do not 
satisfy $\phi( \TT_\qb ) \subseteq \TT_\qb$ since $\phi(X_l^{-1}) = \phi(X_l)^{-1}$ 
belongs to $\TT_\qb $
only in very special cases.
In Section \ref{aut-q-tor} we prove the 
following result:

\bth{i} Let $\TT_\qb$ be a saturated quantum torus of rank $N$
over an arbitrary base field $\KK$. Let $\db \in \ZNp$ be 
a degree vector. For every bi-finite unipotent automorphism
$\phi$ of the completed quantum torus $\wh{\TT}_{\qb, \db}$, 
there exist $N$ elements
\[
u_1, u_2, \ldots, u_N \in Z(\TT_\qb)^{\geq 1}
\] 
such that $\phi(X_l) = (1+u_l) X_l$
for all $1 \leq l \leq N$,
where $Z(\TT_\qb)^{\geq 1} := Z(\TT_\qb) \cap (\wh{\TT}_{\qb, \db})^{\geq 1}$.
\eth 
\thref{i} is a rigidity result for bi-finite 
unipotent automorphisms of completed saturated quantum tori in the sense 
that it asserts that those kinds of automorphisms are only coming 
from the center of the underlying torus.
We derive the Andruskiewitsch--Dumas conjecture
from a rigidity theorem for a type of unipotent automorphisms
of the algebras $\UU^+_q(\g)$. Every strictly dominant integral coweight 
$\la= \sum_{\al \in \Pi} m_\al \vpi_\al\spcheck$ of $\g$
gives rise to a 
connected $\Zp$-grading of $\UU_q^+(\g)$ 
obtained by assigning degree $m_\al=\lcor \la, \al \rcor$ 
to $E_\al$, where $\{ \vpi_\al\spcheck \mid \al \in \Pi \}$ 
are the fundamental coweights of $\g$. For $m \in \Zp$, denote by 
$\UU_q^+(\g)^{\geq m}$ the space of elements  
of degree $\geq m$. We call an automorphisms 
$\Phi$ of $\UU_q^+(\g)$ $\la$-unipotent if it satisfies 
\[
\Phi(E_\al) - E_\al \in \UU_q^+(\g)^{\geq \lcor \la, \al \rcor + 1}, \; \; 
\forall \al \in \Pi.
\]
\bth{i2} Let $\g$ be a simple Lie algebra of rank $r>1$, $\KK$ 
an arbitrary base field, $q$ a deformation parameter that is 
not a root of unity. For every strictly dominant integral coweight $\la$ 
the only $\la$-unipotent automorphism 
of $\UU^-_q(\g)$ is the identity automorphism.
\eth
Our strategy for the proof of this theorem is as follows. 
The algebras $\UU_q^+(\g)$ are examples of a type of iterated Ore 
extensions called Cauchon--Goodearl--Letzter extensions and
the Cauchon procedure 
of deleting derivations \cite{Ca} can be used to embed them 
into quantum tori, see \S \ref{2.2} for details. This is 
not yet sufficient to relate the $\la$-unipotent automorphisms of $\UU_q^+(\g)$ 
to bi-finite unipotent automorphisms of completed quantum tori. 
For this we apply a recent result of Geiger and the author \cite{GeY}
stating that one can change the generators of those quantum tori
so they become quantum minors in $\UU^+_q(\g)$. Recall that a 
quantum affine space algebra is an algebra with generators 
$X_1, \ldots, X_N$ and relations as in \eqref{q-tor}.
The above mentioned result of \cite{GeY} leads to a chain
of embeddings
\begin{equation}
\label{embb}
\AA \subset \UU_q^+(\g) \subset \TT,
\end{equation}
where $\AA$ is a quantum affine space algebra and $\TT$ is the corresponding 
quantum torus (which coincides with the Cauchon quantum torus).
In the case when $\KK$ has characteristic 0 and $q$ is transcendental over $\Qset$
one can also obtain this by applying the  
results of Gei\ss--Leclerc--Schr\"oer \cite{GLS}. 
Using \eqref{embb} we find a relationship between the $\la$-unipotent 
automorphisms of $\UU_q^+(\g)$ and the bi-finite unipotent 
automorphisms of a completion of the (saturated) quantum torus $\TT$. 
Then we apply results from \cite{Cal2,Y4} on the normal elements of the algebras 
$\UU^+_q(\g)$ and a theorem for separation of variables
for these algebras. These results and \thref{i} are used to 
prove that every $\la$-unipotent automorphism $\Phi$ of $\UU_q^+(\g)$ 
satisfies
\[
\Phi(E_\al) = (1 + z_\al) E_\al, \; \; 
\forall \al \in \Pi
\]
for some $z_\al \in Z(\UU_q^+(\g)) \cap \UU_q^+(\g)^{\geq 1}$.
Finally, the structure of the torus invariant height 
one prime ideals of $\UU_q^+(\g)$ from \cite{J3,Go,Y4} is used 
to establish that $z_\al = 0$, $\forall \al \in \Pi$. 
The proof of \thref{i2} is given in Section \ref{UqUnip}.
Section \ref{sAD} contains the proof of 
the Andruskiewitsch--Dumas conjecture. It is based 
on \thref{i2} and an intermediate classification of the family
of automorphisms of $\UU_q^+(\g)$ that map the subspace
$\Span \{ E_\al \mid \al \in \Pi\}$ to itself.

In Section \ref{multi} we prove an extension of 
\cjref{AD}, which classifies the automorphism groups 
of the 2-cocycle twists of the algebras $\UU_q^+(\g)$,
again in full generality.
All proofs in the paper are carried out in such a way 
so they easily extend to the twisted case.
We do not go straight to the twisted case to avoid 
technicalities, which will obscure the main ideas.

In Section \ref{multi} the results on 
automorphism groups are also applied to obtain a full solution 
of the isomorphism problem for the family of all algebras 
obtained by 2-cocycle twists from the algebras $\UU^+_q(\g)$
for simple Lie algebras $\g$. In particular, it is shown that 
\[
\UU_q^+(\g_1) \cong \UU_q^+(\g_2) \quad \Leftrightarrow \quad 
\g_1 \cong \g_2
\]
for all base fields $\KK$, non-root of unity $q$ and 
simple Lie algebras $\g_1$, $\g_2$. The idea to apply 
the results on automorphism groups to this isomorphism problem 
was suggested by Len Scott.

The methods of this paper have a very broad range of applications
to the investigation of automorphism groups of noncommutative algebras.
They provide a procedure to deal with individual automorphisms
or analyze the full automorphism groups of those algebras 
$R$ that satisfy 
\begin{equation}
\label{me}
\AA \subset R \subset \TT
\end{equation}
for some quantum torus $\TT$ and the corresponding quantum affine 
space algebra $\AA$. This is done by using the above mentioned relationship
between the automorphisms of $R$ and the bi-finite unipotent 
automorphisms of a completion of $\TT$, and then applying the 
rigidity from \thref{i}. In its most general form the former relationship
is stated in \cite[Proposition 3.3]{Y8} in connection to one such application. 
There are very large classes of algebras $R$ that satisfy \eqref{me}. 
For example all quantum cluster algebras. (The above procedure 
in this case deals with the full automorphism group, not 
just maps that take clusters to clusters.) In a recent preprint
\cite{GY} K. Goodearl and the author proved that the property \eqref{me} 
is satisfied by all algebras in the large, axiomatically defined 
class of iterated Ore extensions called Cauchon--Goodearl--Letzter 
extensions \cite{Ca,GL}. There are particular families of algebras in the above classes 
for which the automorphism groups have been of interest.
In \cite{Y8} we apply the methods of this 
paper to prove the Launois--Lenagan conjecture \cite{LaL1}
that for all integers $N \geq 2$ the automorphism group of the algebra 
$R_q[M_N]$ of quantum matrices of size $N \times N$ is isomorphic to a semidirect product 
of the torus $\Tset^{2 N-1}$ and a copy of $\Zset_2$ corresponding
to the transpose automorphism. It was proved \cite{ACh,LaL2} for $N=2$ and $3$, 
and was open for all $N >3$. Other particular families of algebras 
to which the procedure is applicable and the automorphism groups 
have been of interest include quantum groups $R_q[G]$, \cite{J1,HLT} and   
the quantum Schubert cell algebras $\UU^+[w]$, \cite{DKP,L,Ja}. 

Although the above procedure makes sense for commutative algebras $R$
(and thus for classical cluster algebras), in those cases \thref{i} 
does not produce sufficient restrictions on the possible 
form of the automorphisms of $R$.
However, there is a Poisson version of \thref{i} about rigidity 
of automorphisms of Poisson tori. This and some of its applications 
will be described in another publication.
\medskip
\\
\noindent
{\bf Acknowledgements.} I am grateful to Len Scott for suggesting 
the application of the main results of the paper to 
the isomorphism problem for the algebras $\UU^+_q(\g)$
(\thref{isom}) and
to Ken Goodearl for 
his very helpful comments on the first version
of the manuscript,
advice on terminology, and for letting me know of the 
paper \cite{Ar}.
I would also like to thank Jaques Alev and 
Nicol\'as Andruskiewitsch for valuable discussions
and the referee for the careful reading of the 
manuscript and his/her suggestions which improved 
the exposition.

The author was supported in part by NSF grant DMS-1001632.
\sectionnew{The algebras $\UU_q^\pm(\g)$}
\lb{qalg}
\subsection{Quantized universal enveloping algebras}
\label{2.1}
We will mostly follow the notation of Jantzen's book \cite{Ja}. 
Assume that $\g$ is a complex simple Lie algebra with Weyl 
group $W$, set of simple roots $\Pi$, and Dynkin diagram $\Ga$. 
We will identify the set of vertices of $\Ga$ 
with $\Pi$. For $\al \in \Pi$ denote by $\vpi_\al$ and 
$s_\al \in W$ the corresponding fundamental weight and simple 
reflection. Let $\QQ$ and $\PP$ be the root and weight lattices of $\g$.
Let $\QQ_+= \Zp \Pi$ and
$\PP_+ = \Zp \{ \vpi_\al \mid \al \in \Pi\}$ be the set 
of dominant integral weights. 
The support of $\mu = \sum_{\al \in \Pi} m_\al \vpi_\al \in \PP$ 
is defined by 
$\Supp (\mu) : = \{ \al \in \Pi \mid m_\al \neq 0 \}$.
Denote by
$\PP_{++}\spcheck = \{ \sum m_\al \vpi_\al\spcheck 
\mid m_\al \in \Zset_+, \forall \al \in \Pi \}$
the set of strictly dominant integral coweights of $\g$,
where $\vpi_\al\spcheck$ are the fundamental coweights of $\g$.
Let $\lcor.,. \rcor$ be the invariant bilinear form 
on $\Rset \Pi$ normalized by $\lcor \al, \al \rcor = 2$ for short roots 
$\al \in \Pi$. 

Throughout the paper $\KK$ will denote a base field
(of arbitrary characteristic) and $q \in \KK^*$ 
an element that is not 
a root of unity. The quantized universal enveloping 
algebra $\UU_q(\g)$ of $\g$ is the $\KK$-algebra 
with generators $\{ K_\al^{\pm 1}, E_\al, F_\al \mid  
\al \in \Pi \}$ and relations \cite[\S 4.3]{Ja}.
Let $\UU_q^+(\g)$ and $\UU_q^-(\g)$
be the subalgebras of $\UU_q(\g)$ generated by 
$\{E_\al\mid \al \in \Pi \}$ and
$\{F_\al\mid \al \in \Pi \}$.
There is a unique automorphism $\omega$ 
of $\UU_q(\g)$ such that 
\begin{equation}
\label{om}
\om(E_\al) = F_\al, \; 
\om(F_\al) = E_\al, \; 
\om(K_\al) = K_{\al}^{-1}, \; \; 
\forall \, \al \in \Pi.
\end{equation}
It restricts to 
an isomorphism $\om \colon \UU_q^\pm(\g) \to \UU_q^\mp(\g)$.
We will work with the algebra $\UU^-_q(\g)$ since we use 
results from \cite{GeY,Y6} which will need an appropriate 
reformulation for $\UU^+_q(\g)$.
The algebra $\UU_q^-(\g)$ is the $\KK$-algebra 
with generators $\{F_\al\mid \al \in \Pi \}$ and 
the following quantum Serre relations:
\begin{equation}
\sum_{j=0}^{1-a_{\al \al'}} (-1)^j
\begin{bmatrix} 
1-a_{\al \al'} \\ j
\end{bmatrix}_{q_\al}
      (F_\al)^j F_{\al'} (F_\al)^{1-a_{\al \al'}-j} = 0, 
\; \; \forall \al \neq \al' \in \Pi,
\label{S1}
\end{equation}
where 
\begin{equation}
\label{a}
a_{\al \al'} = 2 \lcor \al, \al' \rcor/ \lcor \al, \al \rcor
\end{equation}
and $q_\al = q^{\lcor \al, \al \rcor/2}$. Here 
$[m]_q = (q^m - q^{-m})/(q-q^{-1})$ for $m \geq 1$, $[0]_q =1$, and
$[m]_q! = [0]_q \ldots [m]_q$, 
$\begin{bmatrix} m \\ j \end{bmatrix}_q = [m]_q!/ [j]_q! [m-j]_q!$
for $j\leq m \in \Zp$. Denote 
\begin{equation}
\label{numbers}
r := \rank(\g), \; N:= (\dim \g - r)/2, \;  
\mbox{and} \; 
\Tset^r := (\KK^*)^r.
\end{equation}
In other words $N$ is the number of positive roots of $\g$.   
Then $\GKdim \UU_q^\pm(\g) = N$. (This follows for instance 
from the iterated Ore extension presentation \eqref{Umin-iter} below
and the isomorphism $\UU^+_q(\g) \cong \UU^-_q(\g)$.) 
The algebra $\UU_q(\g)$ is $\QQ$-graded by assigning 
$E_\al$, $F_\al$, and $K_\al^{\pm 1}$ weights $\al$, $-\al$, and $0$. 
For $\ga \in \QQ$, the corresponding graded component of 
$\UU_q(\g)$ will be denoted by $\UU_q(\g)_\ga$. The grading 
gives rise to the following $\Tset^r$-action on $\UU_q(\g)$ by
algebra automorphisms:
\begin{equation}
\label{torus-act}
t \cdot u = t^\ga u, \quad u \in \UU_q(\g)_\ga, \ga \in \QQ
\end{equation}
in terms of the characters
\[
t \mt 
t^\ga := \prod_{ \al \in  \Pi} t_\al^{\lcor \ga, \vpi_\al \rcor}, \quad 
t =(t_\al)_{\al \in \Pi} \in \Tset^r.
\] 
\subsection{Cauchon's procedure of deleting derivations and $\UU_q^\pm(\g)$}
\label{2.2}
Consider an iterated Ore extension 
\begin{equation} 
\label{itOre}
R := \KK[x_1][x_2; \sig_2, \delta_2] \ldots [x_N; \sig_N, \delta_N],
\end{equation}
where for $l \in [2,N]$, $\sig_l$ is an automorphism 
and $\delta_l$ is a (left) $\sig_l$-skew derivation of 
$R_{l-1}:=\KK[x_1][x_2; \sig_2, \delta_2] 
\ldots [x_{l-1}; \sig_{l-1}, \delta_{l-1}]$.
Here and below for $m \leq n \in \Zset$, we set $[m,n] =\{m, \ldots, n \}$.

\bde{CGL} An iterated Ore extension $R$ given by \eqref{itOre} 
is called a Cauchon--Goodearl--Letzter (CGL) extension if it is 
equipped with an action of the torus $\Tset^r := (\KK^*)^r$
by algebra automorphisms satisfying the following conditions:

(i) For all $1 \leq k < l \leq N$, $\sig_l(x_k) = q_{lk} x_k$
for some $q_{lk} \in \KK^*$.

(ii) For every $l \in [2,N]$, $\delta_l$ is a locally nilpotent 
$\sig_l$-skew derivation of $R_{l-1}$. 

(iii) The elements $x_1, \ldots, x_N$ are $\Tset^r$-eigenvectors 
and the set $\{ c \in \KK \mid \exists t \in \Tset^r, \; 
t \cdot x_1 = c x_1 \}$ is infinite.

(iv) For every $l \in [2,N]$ there exists $t_l \in \Tset^r$ such that 
$t_l \cdot x_l = q_l x_l$ for some $q_l \in \KK^*$ which is 
not a root of unity, and $t_l \cdot x_k = q_{lk} x_k$, 
$\forall k \in [1,l-1]$ (i.e., $\sig_l = (t_l \cdot)$ as elements 
of $\Aut(R_{l-1})$, $\forall l \in [2,N]$). 
\ede 

We note that for all CGL extensions, 
$\sig_l \delta_l =q_l \delta_l \sig_l$, $\forall l \in [2,N]$.
In this setting Cauchon \cite{Ca} iteratively constructed
$N$-tuples of nonzero elements 
\[
(x^{(m)}_1, \ldots, x^{(m)}_N)
\] 
of the division ring of fractions $\Fract(R)$ of $R$ for
$m= N+1, \ldots, 1$. First,
\[
(x^{(N+1)}_1, \ldots, x^{(N+1)}_N) := 
(x_1, \ldots, x_N).
\]
The other $N$-tuples are obtained recursively from the formula
\begin{equation}
\label{new-x}
x^{(m)}_j := 
\begin{cases}
x^{(m+1)}_j, 
& \mbox{if} \; \; j \geq m 
\\
\sum_{n=0}^\infty \frac{(1- q_m)^{-n}}{(n)_{q_m}!} 
\Big[ \delta_m^n \sig^{-n}_m \left(x^{(m+1)}_j \right) \Big]
\left(x^{(m+1)}_m \right)^{-n}, 
& \mbox{if} \; \; j < m
\end{cases}
\end{equation}
for $m =N, \ldots, 2$. Here $(n)_q! = (0)_q \ldots (n)_q$,
$(0)_q=1$, and $(n)_q= (1-q^n)/(1-q)$ for $n\geq 1$.
This process is called {\em{Cauchon's 
procedure of deleting derivations}}. The terminology comes
from the following fact proved by Cauchon \cite{Ca}: the subalgebra of 
$\Fract(R)$ generated by the $m$-th $N$-tuple of elements 
is isomorphic to an iterated Ore extension of the form 
\eqref{itOre}, where the derivations $\de_m, \ldots, \de_N$
are no longer present. Denote the final $N$-tuple of elements
$(\ol{x}_1, \ldots, \ol{x}_N):= (x^{(2)}_1, \ldots, x^{(2)}_N)$
and the subalgebra of $\Fract(R)$ generated those elements and 
their inverses
\begin{equation}
\label{C-tor} 
\TT = \KK \lcor \ol{x}_1^{\, \pm 1}, \ldots, \ol{x}_N^{\, \pm 1} \rcor.
\end{equation}
Let $q_{ll}= 1$ for $l \in [1,N]$, $q_{kl} = q^{-1}_{lk}$ 
for $1 \leq k < l \leq N$, and $\qb := (q_{kl})_{k,l=1}^N \in M_N(\KK^*)$.
Cauchon proved \cite{Ca} that 
\begin{equation}
\label{Cau}
R \subset \TT \; \; 
\mbox{and the map} \; \; 
\eta \colon \TT_{\qb} \to \TT \; \; 
\mbox{given by} \; \; \eta(X_l) = \ol{x}_l 
\; \; \mbox{is an isomorphism},
\end{equation}
recall \eqref{q-tor}. The CGL extension $R$ is called {\em{torsion-free}} 
if the subgroup of $\KK^*$ generated by $q_{kl}$ for $1 \leq k < l \leq N$ 
is torsion-free. In such a case $\TT_\qb$ is a saturated quantum 
torus as noted in the introduction.

Denote the braid group of $\g$ by $\B_\g$ and its standard 
set of generators by $\{T_\al \mid \al \in \Pi\}$.
Let $w_0$ be the longest element of the Weyl group $W$ of $\g$.
A word $\ib = (\al_1, \ldots, \al_N)$ in the alphabet 
$\Pi$ is called a reduced word for $w_0$ if $s_{\al_1} \ldots s_{\al_N}$ 
is a reduced expression of $w_0$, recall \eqref{numbers}. 
For such a reduced word $\ib$ define
\begin{equation}
\label{leq1}
w_0(\ib)_{\leq l} := s_{\al_1} \ldots s_{\al_l}, \; \; l \in [0,N]
\end{equation}
and the Lusztig root vectors
\begin{equation}
F_{\be_l} := 
T_{\al_1} \ldots T_{\al_{l-1}} 
(F_{\al_l}), \; \;
\mbox{where} \; \; 
\beta_l := w_0(\ib)_{\leq (l-1)} \al_{i_l}, \; \; 
l \in [1, N],
\label{rootv}
\end{equation}
see \cite[\S 39.3]{L}.
Here we use Lusztig's action 
of $\B_\g$ on $\UU_q(\g)$ 
in the version given in \cite[\S 8.14]{Ja} by Eqs. 
8.14 (2), (3), (7), and (8). We will need the following 
property (see \cite[Proposition 8.20]{Ja}):
\begin{equation}
\label{eqa}
\mbox{if} \; \; \be_l = \al \in \Pi \; \;
\mbox{for some} \; \; l \in [1,N], \; \; \mbox{then}
\; \;  
F_{\be_l} = F_\al.
\end{equation}

Given an algebra $B$, a subalgebra 
$B'$ of $B$, $x \in B$, an automorphism $\sig$ of $B'$, and a 
(left) $\sig$-derivation $\delta$ of $B'$, we will say that 
\[
B = B'[x;\sig, \delta]
\] 
is an Ore extension presentation of $B$ if the map 
$\psi \colon B'[y; \sig, \delta] \to B$ 
given by $\psi(b') = b'$, $b' \in B'$ and 
$\psi(y) = x$ is an algebra isomorphism. 
The Levendorskii--Soibelman 
straightening law 
\begin{multline}
\label{LS}
F_{\be_l} F_{\be_k} - 
q^{ - \lcor \be_l, \be_k \rcor }
F_{\be_k} F_{\be_l}  \\
= \sum_{ {\bf{m}} = (m_{k+1}, \ldots, m_{l-1}) \in \Zp^{l-k-2} }
c_{\bf{m}} (F_{\be_{l-1}})^{m_{l-1}} \ldots (F_{\be_{k+1}})^{m_{k+1}},
\; \; c_{\bf{m}} \in \KK,
\end{multline}
for $1 \leq k < l \leq N$ (see e.g. \cite[Proposition I.6.10]{BG}) 
is used to associate to each reduced word $\ib$ for $w_0$ 
an iterated Ore extension presentation of $\UU_q^-(\g)$. For $l \in [1,N]$ 
choose an element $t_l \in \Tset^r$ such that 
$t_l^{\be_k} = q^{\lcor \be_k, \be_l \rcor }$ 
for all $k \in [1,l]$, cf. \eqref{torus-act}.
Let $\UU^-[w_0(\ib)_{\leq l}]$ be the subalgebra 
of $\UU_q^-(\g)$ generated by $F_{\be_1}, \ldots, F_{\be_l}$
for $l \in [0,N]$. 
(This is nothing but the quantum Schubert cell algebra
of De Concini--Kac--Procesi \cite{DKP} and Lusztig \cite{L}
associated to $w_0(\ib)_{\leq l}$.) 
Then $\UU^-[w_0(\ib)_{\leq 0}] = \KK$, $\UU^-[w_0(\ib)_{\leq N}] = \UU_q^-(\g)$,
and for all $l \in [1,N]$ we have the Ore extension presentations 
\[
\UU^-[w_0(\ib)_{\leq l}] = \UU^-[w_0(\ib)_{\leq (l-1)}][F_{\be_{l-1}}, \sig_l, \delta_l].
\]
Here 
\begin{equation}
\label{sig}
\sig_l := (t_l \cdot) \in \Aut( \UU^-[w_0(\ib)_{\leq (l-1)}]) 
\end{equation}
for
the element $t_l \in \Tset^r$ constructed above and the restriction 
of the action \eqref{torus-act} to $\UU^-[w_0(\ib)_{\leq (l-1)}]$. 
The skew derivation $\de_l$ is given by 
\begin{equation}
\label{delta}
\delta_l(x) := F_{\be_l} x - q^{\lcor \be_l, \ga \rcor} x F_{\be_l}, 
\; \; 
x \in \UU^-[w_0(\ib)_{\leq (l-1)}]_\ga, \ga \in \QQ,
\end{equation}
recall \eqref{LS}.
By composing those presentations, one associates
to each reduced word $\ib$ for $w_0$ the
iterated Ore extension presentation of $\UU_q^-(\g)$
\begin{equation}
\label{Umin-iter}
\UU_q^-(\g) = \KK [F_{\be_1}] [F_{\be_2}; \sig_2, \delta_2] \ldots 
[F_{\be_N}; \sig_N, \delta_N ].
\end{equation}
This is a torsion-free CGL extension for
the following choice of the coefficients $q_{lk}$, $q_l$:
\begin{equation}
\label{coeff}
q_{lk} = q^{- \lcor \be_l, \be_k \rcor}, 
\; \; 1 \leq k < l \leq N, \quad
q_l = q_{\al_l}^{-2}, \; \; l \in [1,N],
\end{equation}
see \cite{MC}.
\subsection{Separation of variables for $\UU_q^\pm(\g)$ and height one primes}
\label{2.3}
Recall that a $\UU_q(\g)$-module $V$ is called a type one module if it equals 
the sum of its $q$-weight spaces defined by 
\[
V_\mu := \{ v \in V \mid K_\al v = q^{ \lcor \mu, \al \rcor} v, \; \; 
\forall \, \al \in \Pi \}, \; \mu \in \PP.
\]
The irreducible finite dimensional 
type one $\UU_q(\g)$-modules are parametrized by 
the dominant integral weights of $\g$,
see \cite[Theorem 5.10]{Ja}.
Denote by $V(\la)$ 
the irreducible $\UU_q(\g)$-module with highest weight 
$\la \in \PP_+$, and fix a highest weight vector $v_\la$ 
of $V(\la)$.
The braid group $\B_\g$ acts 
on $V(\la)$ by \cite[Eq. 8.6 (2)]{Ja}. This action 
is compatible with the one on $\UU_q(\g)$ and 
in particular satisfies 
$T_w (V(\la)_\mu) = V(\la)_{w \mu}$, 
$\forall w \in W$, $\la \in \PP_+$, $\mu \in \PP$.
Because of this, for all $w \in W$ there exists a unique element 
$\xi_{w,\la} \in (V(\la)^*)_{- w \la}$ 
such that $\xi_{w,\la} (T^{-1}_{w^{-1}} v_\la) = 1$,
where dual modules are formed using the antipode of 
the Hopf algebra structure on $\UU_q(\g)$. 
Recall that $w_0$ denotes the longest element of $W$.
Given $w \in W$,
denote the matrix coefficients
\begin{equation} 
\label{c-notation}
e^\la_w \in (\UU_q(\g))^*, \quad 
e^\la_w(x) := \xi_{w,\la}( x T^{-1}_{w_0^{-1}} v_\la ), \; \; x \in \UU_q(\g)
\end{equation}
called quantum minors in the case when $\la$ is a fundamental weight.
We will need their counterparts in $\UU_q^-(\g)$. Given 
$\ga \in \QQ_+ \backslash \{ 0 \}$, denote
$n(\ga) = \dim \UU_q^+(\g)_\ga= \dim \UU_q^-(\g)_{-\ga}$
and fix a pair of dual bases 
$\{u_{\ga, j} \}_{j=1}^{n(\ga)}$ and 
$\{u_{-\ga, j} \}_{j=1}^{n(\ga)}$
of $\UU_q^+(\g)_\ga$ and $\UU_q^-(\g)_{-\ga}$ 
with respect to the Rosso--Tanisaki form, see \cite[Ch. 6]{Ja}. 
The universal $R$-matrix corresponding of $\UU_q(\g)$ (without its 
semisimple part) is given by
\begin{equation}
\label{RR}
\RR := 1 \otimes 1 + \sum_{\ga \in \QQ_+, \ga \neq 0} \sum_{j=1}^{n(\ga)} 
u_{\ga, j} \otimes u_{- \ga, j} \in \UU_q^+(\g) \wh{\otimes} \UU_q^-(\g).
\end{equation}
Here $\UU_q^+(\g) \wh{\otimes} \UU_q^-(\g)$ denotes 
the completion of $\UU_q^+(\g) \otimes \UU_q^-(\g)$ 
with respect to the filtration \cite[\S 4.1.1]{L}.
There is a unique graded algebra antiautomorphism 
$\tau$ of $\UU_q(\g)$ given by 
\begin{equation}
\label{tau}
\tau(E_\al) = E_\al,
\,
\tau(F_\al) = F_\al, 
\, 
\tau(K_\al) = K_\al^{-1}, \; \; 
\forall \, \al \in \Pi,
\end{equation}
see \cite[Lemma 4.6(b)]{Ja}. The counterparts of $e_w^\la$ 
in $\UU_q^-(\g)$ are the elements 
\begin{equation}
\label{bel}
b_w^\la := (e^\la_{w,w_0} \tau \otimes \id) \RR^w \in \UU_q^-(\g)_{-(w-w_0)\la}, 
\; \; \la \in \PP_+.
\end{equation}
They play a key role in the description of the 
spectra of the algebras $\UU_q^-(\g)$, see 
\cite[Theorem 3.1]{Y6}.
A more conceptual way to define them is via a family 
of homomorphisms which realize the quantum Schubert cell 
algebras as quotients of quantum function algebras, 
see \cite[Theorem 3.6]{Y2}.

The elements $b_1^\la$, $\la \in \PP_+$ 
are normal elements of $\UU_q^-(\g)$
\begin{equation}
\label{norm}
b_1^\la u
= 
q^{ \lcor (1+w_0)\la, \ga \rcor }
u b_1^\la, \; \; 
\forall \, 
u \in \UU_q^-(\g)_{-\ga},
\ga \in \QQ_+, 
\end{equation}
and satisfy 
\begin{equation}
\label{bprod}
b_1^{\la} b_1^{\la'} = b_1^{\la'} b_1^\la = 
q^{\lcor \la, (1-w_0)\la' \rcor } b_1^{\la + \la'}
\end{equation}
for $\la, \la' \in \PP_+$, see \cite[Eqs. (3.30)-(3.31)]{Y4}.
(Note that for all $\la, \la' \in \PP$, 
$\lcor \la, (1-w_0)\la' \rcor =
\lcor \la', (1-w_0)\la \rcor$.)
Denote by $\NN_q^-(\g)$ the subalgebra 
of $\UU_q^-(\g)$ generated (and hence spanned) by 
$b^\la_1$, $\la \in \PP_+$. 

Finally, for a reduced word $\ib$ of $w_0$ 
define the following subset of $\Zp^{N}$:
\begin{equation}
\label{Hset}
\HH_\ib :=
\{ (j_1, \ldots, j_N) \in \Zp^N 
\mid \forall \al \in \Pi, \; \exists k \in [1,N] 
\; \; \mbox{such that} \; \; 
\al_k = \al \; \; \mbox{and} \; \; 
j_k = 0 \}.
\end{equation}
In other words, for every simple root $\al$ we consider the set $\{ k \in [1,N] \mid \al_k = \al\}$ 
and require that there exists an index $k$ in this set such that $j_k =0$.

We will need the following results from \cite{Y4} describing the structure 
of the algebras $\UU^-_q(\g)$ as $\NN^-_q(\g)$-modules and the set of height 
one $\Tset^r$-prime ideals of $\UU_q^-(\g)$.

\bth{1a} 
\cite[Theorems 5.1, 5.4 and 6.19, and Proposition 6.9]{Y4} 
For all simple Lie algebras $\g$, base fields $\KK$, and 
$q \in \KK^*$ not a root of unity, the following hold:

(i) The height one $\Tset^r$-invariant prime ideals 
of $\UU_q^-(\g)$ are precisely the ideals 
\[
\UU_q^-(\g) b_1^{\vpi_\al}, \; \; \al \in \Pi.
\]

(ii) All normal elements of $\UU_q^-(\g)$ 
belong to $\NN_q^-(\g)$. 
The algebra $\NN_q^-(\g)$ is a polynomial algebra 
in the generators $\{ b^{\vpi_\al}_1 \mid \al \in \Pi\}$,
i.e., $\{ b^\la_1 \mid \la \in \PP_+\}$ is a $\KK$-basis 
of $\NN^-_q(\g)$.

(iii) The algebra 
$\UU_q^-(\g)$ is a free left (and right) 
$\NN_q^-(\g)$-module with basis
\begin{equation}
\label{basis}
\{ F_{\be_N}^{j_N} \ldots F_{\be_1}^{j_1} 
\mid (j_1, \ldots, j_N) \in \HH_\ib
\}. 
\end{equation}
\eth
When $\charr \KK=0$ and $q$ is transcendental over $\Qset$, part (i) 
of the theorem follows from results of Gorelik \cite{Go} 
and Joseph \cite{J3}. Under those conditions on $\KK$ and $q$, 
part (ii) of the theorem was proved by Caldero in \cite{Cal2}. 
The second part of the theorem establishes that $\NN^-_q(\g)$ 
is precisely the subalgebra of $\UU_q^-(\g)$ generated 
by all normal elements of $\UU_q^-(\g)$. The third part of the 
theorem can be viewed as a result for separation variables
for the algebras $\UU_q^-(\g)$, see \cite[Section 5]{Y4}.
\subsection{Cauchon's procedure for $\UU_q^\pm(\g)$ and quantum minors}
\label{2.4}
Given a reduced word $\ib = (\al_1, \ldots, \al_N)$ for $w_0$, 
consider the CGL extension presentation \eqref{Umin-iter} 
of $\UU_q^-(\g)$. Denote by $(\ol{F}_{\ib,1}, \ldots, \ol{F}_{\ib,N})$
the final $N$-tuple $(\ol{x}_1, \ldots, \ol{x}_N)$ 
of the Cauchon procedure of deleting derivations applied to it. Define
a successor function $s \colon [1,N] \to [1,N] \sqcup \{ \infty \}$
associated to $\ib$ by
\begin{equation}
\label{succ}
s(l) = \min \{ k \mid k > l, \al_k = \al_l \}, \; 
\mbox{if} \; \exists k > l \; \; \mbox{such that} 
\; \al_k = \al_l, \quad s(l) = \infty, \; 
\mbox{otherwise}.
\end{equation}
Define the quantum minors 
\begin{equation}
\label{minors}
\De_{\ib, l} := b^{\vpi_{\al_l}}_{ w_0(\ib)_{\leq (l -1)}},
\; \; l \in [1,N].
\end{equation}
The following result from \cite{GeY} expressing the 
Cauchon elements $\ol{F}_{\ib, l}$ in terms of the 
quantum minors $\De_{\ib,l}$ will be needed later:

\bth{1b} \cite[Theorem 3.1]{GeY} For all 
simple Lie algebras $\g$, base fields $\KK$,
$q \in \KK^*$ not a root of unity, and 
reduced words $\ib$ for $w_0$, the elements 
$\ol{F}_{\ib,1}, \ldots, \ol{F}_{\ib, N} \in \Fract(\UU_q^-(\g))$ 
from the Cauchon deleting derivation procedure for the torsion-free
CGL presentation \eqref{Umin-iter} of $\UU_q^-(\g)$ are given by
\[
\ol{F}_{\ib, l} = 
\begin{cases}
(q_{\al_l}^{-1} - q_{\al_l})^{-1} 
\De_{\ib, s(l)}^{-1} \De_{\ib, l} , 
&\mbox{if} \; \; s(l) \neq \infty 
\\
(q_{\al_l}^{-1} - q_{\al_l})^{-1} \De_{\ib, l}, 
&\mbox{if} \; \; s(l) = \infty.
\end{cases}
\]
\eth
\sectionnew{Automorphisms of completed quantum tori} 
\label{aut-q-tor}
\subsection{Quantum tori}
\label{3.1} 
Let $\KK$ be an arbitrary field
and $\qb = (q_{kl})_{k,l=1}^N \in M_N(\KK^*)$ a multiplicatively
skew-symmetric matrix. Recall the definition \eqref{q-tor} of the 
quantum torus $\TT_\qb$. Denote 
\[
X^{(j_1, \ldots, j_N)}:= X_1^{j_1} \ldots X_N^{j_N}
\in \TT_\qb, 
\; \; (j_1, \ldots, j_N) \in \ZN. 
\]
Let $\{ e_1, \ldots, e_N \}$ be the standard basis of $\ZN$.
Thus 
\begin{equation}
\label{e-gen}
X^{e_k} = X_k, \; \; k \in [1,N].
\end{equation}
The quantum torus $\TT_\qb$ has the $\KK$-basis
$\{ X^f \mid f \in \ZN \}$. Recall the definition 
\eqref{Kerq} of the multiplicative kernel of $\qb$. 
We have
\begin{equation}
\label{Cent}
Z(\TT_\qb) = \Span \{ X^f \mid f \in \Ker(\qb) \}.
\end{equation}
It is straightforward to show that each of the following 
three conditions is equivalent to $\TT_\qb$ 
being saturated, recall \eqref{torsion-free}:
\begin{gather}
\Ker(\qb) = (\Ker(\qb) \otimes_\Zset \Rset) \cap \ZN,
\label{c1}
\\
\mbox{for} \; \; f \in \ZN, n \in \Zset_+, \; \; 
X^{nf} \in Z(\TT_\qb) \Rightarrow X^f \in Z(\TT_\qb),
\label{c2}
\\
\mbox{for} \; \; u \in \TT_\qb, n \in \Zset_+, \; \; 
u^n \in Z(\TT_\qb) \Rightarrow u \in Z(\TT_\qb). 
\label{c3}
\end{gather}
It follows from either of the two conditions \eqref{c2} and \eqref{c3}
that the property of a torus $\TT_\qb$ being 
saturated is independent of the choice of 
generators. 

%

An $N$-tuple $\db = (d_1, \ldots, d_N) \in \Zp^N$ 
will be called 
a {\em{degree vector}}. Such will be used to define 
a completion of $\TT_\qb$ as follows. Define the homomorphism 
\[
D \colon \ZN \to \Zset, \; \;  
D(j_1, \ldots, j_N) = d_1 j_1 + \ldots + d_N j_N
\] 
and the $\Zset$-grading
\[
\TT_\qb = \oplus_{ m \in \Zset} \TT_\qb^m, \; \; 
\TT_\qb^m = \Span \{ X^f \mid D(f) = m \}.
\] 
Consider the associated valuation $\vv \colon \TT_\qb \to \Zset \sqcup \{ \infty \}$
given by
\[
\vv( u_m + u_{m+1} + \ldots + u_{m'}) := m \; \; 
\mbox{for} \; \; 
u_j \in \TT_\qb^j, j \in [m,m'], u_m \neq 0.
\] 
The completion of $\TT_\qb$ with respect to this valuation is given by 
\[
\wh{\TT}_{\qb, \db} : = 
\{ u_m + u_{m+1} + \ldots \mid m \in \Zset, 
u_j \in \TT_\qb^j \; \mbox{for} \; j \geq m \}. 
\]
For $m \in \Zset$, denote 
\[
\wh{\TT}_{\qb, \db}^{\geq m} := 
\{ u_m + u_{m+1} + \ldots \mid u_j \in \TT_\qb^j \; \mbox{for} \; j \geq m \}
= \{ u \in \wh{\TT}_{\qb, \db} \mid \vv(u) \geq m \}.
\]
For every graded subalgebra $R$ of $\wh{\TT}_{\qb,\db}$,
set $R^{\geq m}= R \cap \wh{\TT}_{\qb,\db}^{\geq m}$.
It is straightforward to verify that the group of units of $\wh{\TT}_{\qb, \db}$
is given by
\begin{equation}
\label{unit1}
U(\wh{\TT}_{\qb, \db}) = \{ c X^{f} + u \mid 
c \in \KK^*, f \in \ZN, u \in \wh{\TT}_{\qb, \db}^{\geq D(f) + 1} \},
\end{equation}
and that for an element $c X^{f} + u $ as in the right hand side 
of \eqref{unit1} we have
\begin{equation}
\label{unit2}
(c X^{f} + u )^{-1} = c^{-1} X^{-f} ( 1 +  c^{-1} u X^{-f})^{-1} =
\sum_{m=0}^\infty (-1)^m c^{-m -1} X^{-f} (u X^{-f})^m.
\end{equation}

\bde{unip} A continuous automorphism $\phi$ of $\wh{\TT}_{\qb, \db}$ 
will be called {\em{unipotent}} if 
\[
\phi(X_k) - X_k \in \wh{\TT}_{\qb,\db}^{\geq d_k + 1} \; \;
\mbox{for all} \; k \in [1,N].
\]
A unipotent automorphism $\phi$ of $\wh{\TT}_{\qb, \db}$ 
will be called {\em{finite}} if 
\[
\phi(X_k) \in \TT_\qb \; \; 
\mbox{for all} \; k \in [1,N].
\]
A unipotent automorphism $\phi$ of $\wh{\TT}_{\qb,\db}$ 
will be called {\em{bi-finite}} if both $\phi$ and 
$\phi^{-1}$ are finite.
\ede

Not all finite unipotent automorphisms are bi-finite.
In the single parameter case certain finite unipotent automorphisms 
that are not bi-finite
play an important role in the Berenstein--Zelevinsky work on 
quantum cluster algebras, see \cite[Proposition 4.2]{BZ}. 

\bre{unip-aut}
(i) An automorphism $\phi$ of $\wh{\TT}_{\qb,\db}$ is 
{\em{unipotent}} if and only if 
\begin{equation}
\label{un-2}
\phi(u) - u \in \wh{\TT}_{\qb, \db}^{\geq m+ 1}, \; \; 
\forall u \in \wh{\TT}_{\qb, \db}^{\geq m}, m \in \Zset,
\end{equation}
which is also equivalent to
\[
\vv \big(  (\phi -\id) u \big) \geq \vv(u) + 1, \; \; 
\forall u \in \wh{\TT}_{\qb, \db}.
\]
An endomorphism $\phi$ of $\wh{\TT}_{\qb, \db}$ satisfying 
\eqref{un-2} is a continuous automorphism of $\wh{\TT}_{\qb, \db}$.
The set of all unipotent automorphisms of $\wh{\TT}_{\qb, \db}$
is a subgroup of the group of all continuous automorphisms of $\wh{\TT}_{\qb, \db}$.

(ii) There is an obvious isomorphism between the group of all 
automorphisms $\phi$ of $\Fract(\TT_\qb)$ such that
\begin{equation}
\label{fra}
\phi(X_k) - X_k, \phi^{-1}(X_k) -X_k \in \TT_\qb^{\geq d_k +1}, 
\; \; \forall k \in [1,N] 
\end{equation}
and the group of bi-finite unipotent automorphisms of $\wh{\TT}_{\qb, \db}$.
\ere

\bex{fin} Denote by $\AA_\qb$ the quantum affine space subalgebra 
of $\TT_\qb$ generated by $X_1, \ldots, X_N$. In special 
cases the automorphism groups of such algebras were studied by Alev and Chamarie \cite{ACh}.
Every automorphism $\psi$ of $\AA_q$ such that 
\begin{equation}
\label{q-aff}
\psi(X_k) - X_k \in \AA_\qb^{\geq d_k +1}, \; \; \forall k \in [1,N]
\end{equation}
uniquely extends to a bi-finite unipotent automorphism of 
$\wh{\TT}_{\qb, \db}$. In particular, this applies
to the automorphisms of polynomial algebras 
(the case $q_{kl}=1$ for all $k,l$) satisfying \eqref{q-aff} 
for $\db=(1, \ldots, 1)$.
Such automorphisms appear in various contexts.
\eex

\ble{unip} For all multiplicatively skew-symmetric matrices 
$\qb= (q_{kl})_{k,l=1}^N$ and degree vectors $\db$, 
the set of unipotent automorphisms of 
the completed quantum torus $\wh{\TT}_{\qb, \db}$ is in bijection with the $N$-tuples
$(u_1, \ldots, u_N)$ of elements of $\wh{\TT}_{\qb, \db}^{\geq 1}$
such that 
\[
(1+u_k) X_k (1+u_l) X_l 
= q_{kl}
(1 + u_l) X_l (1 + u_k) X_k 
\]
for all $1 \leq k <l \leq N$.
\ele
\begin{proof} If $\phi$ is a unipotent automorphism 
of $\wh{\TT}_{\qb, \db}$, then the $N$-tuple 
\[
(\phi(X_1)X_1^{-1} -1 , \ldots, \phi(X_N)X_N^{-1} - 1) 
\]
satisfies the required property.

For an $N$-tuple $(u_1, \ldots, u_N)$ with that property, 
first define 
\begin{equation}
\label{ph-inv}
\phi(X_k):= (1+ u_k) X_k \; \; 
\mbox{and} \; \; 
\phi( X_k^{-1} ):= X_k^{-1} (1+u_k)^{-1} \in \wh{\TT}_{\qb, \db}, 
\; \; k \in [1,N],
\end{equation}
cf. \eqref{unit2}. Then extend $\phi$ to $\wh{\TT}_{\qb, \db}$
by multiplicativity
\[
\phi( X^{ j_1 e_1 + \ldots + j_N e_N} ):=
\phi(X_1^{\sign(j_1)})^{|j_1|} \ldots \phi(X_N^{\sign(j_N)})^{|j_N|}, \; \; 
\forall j_1, \ldots, j_N \in \Zset,
\]
linearity, and continuity.
It is straightforward to show that the map $\phi$, 
constructed in this way, is an 
endomorphism of $\wh{\TT}_{\qb, \db}$, which satisfies \eqref{un-2}. 
Thus it is a unipotent automorphism of $\wh{\TT}_{\qb, \db}$, see 
\reref{unip-aut} (i).  
\end{proof}
\bre{2} The algebra of all continuous endomorphisms of $\wh{\TT}_{\qb,\db}$ 
is closely related to the group of the unipotent automorphisms of $\wh{\TT}_{\qb,\db}$. 
Every continuous endomorphism $\psi$ of $\wh{\TT}_{\qb, \db}$ satisfies
\begin{equation}
\label{psi}
\psi(X_k) = c_k (1 + u_k) X^{f_k}, \; \; \forall k \in [1,N]
\end{equation}
for some $f_k \in \ZN$, $c_k \in \KK^*$, $u_k \in \wh{\TT}_{\qb, \db}^{\geq 1}$ 
which have the properties
\begin{align*}
&(1+u_k) X^{f_k} (1+u_l) X^{f_l} 
= q_{kl} (1 + u_l) X^{f_l} (1 + u_k) X^{f_k}, \; \;  \forall k <l \; \; \mbox{and}
\\ 
&\exists b \in \Qset_+  \; \; \mbox{such that} \; \; D(f_k) = b d_k, \; \; 
\forall k \in [1,N].
\end{align*}
(This is easily proved using the continuity of $\psi$ and the fact that all $\psi(X_k)$ 
should be units of $\wh{\TT}_{\qb, \db}$, cf. \eqref{unit1}.)
Furthermore, every $\psi$ as in \eqref{psi}
with $c_k, f_k, u_k$ having the stated properties,
uniquely extends to a continuous endomorphism of $\wh{\TT}_{\qb,\db}$. 

\leref{unip} implies that among all continuous endomorphisms of $\wh{\TT}_{\qb,\db}$,
the unipotent automorphisms of $\wh{\TT}_{\qb, \db}$ are precisely the maps \eqref{psi} 
for which $c_k=1$ and $f_k =e_k$ for all $k \in [1,N]$.
\ere

The next theorem contains the main result in this section.

\bth{m1} Assume that $\KK$ is an arbitrary base field,
$\qb \in M_N(\KK^*)$ a multiplicatively 
skew-symmetric matrix for which the quantum torus
$\TT_\qb$ is saturated, and $\db \in \ZNp$ a degree vector. 
Then for every bi-finite unipotent automorphism
$\phi$ of the completed quantum torus $\wh{\TT}_{\qb, \db}$, 
there exists an $N$-tuple
\[
(u_1, u_2, \ldots, u_N) \; \; \mbox{of elements of } \; \; 
Z(\TT_\qb)^{\geq 1}
\] 
such that $\phi(X_k) = (1+u_k) X_k$
for all $k \in [1,N]$.
\eth 

In light of \reref{unip-aut} (ii) this theorem can be restated
as follows:
\medskip

{\em{For a saturated quantum torus $\TT_\qb$ and $\db \in \Zset_+^N$,
every automorphism $\phi$ of $\Fract (\TT_\qb)$ that 
satisfies \eqref{fra} is given by $\phi(X_k) = (1+u_k) X_k$ for 
some elements $u_1, u_2, \ldots, u_N \in  
Z(\TT_\qb)^{\geq 1}$.}}
\medskip

\thref{m1} is proved in \S \ref{3.3}. 
The statement of the theorem does not hold if one replaces bi-finite with 
finite unipotent automorphisms of completed saturated quantum tori. 
Counterexamples are provided by the Berenstein--Zelevinsky automorphisms 
in \cite[Proposition 4.2]{BZ}. Recalling \exref{fin} we note the following:

\bco{fin2} Under the above assumptions on $\qb$ and $\db$, for 
every automorphism $\psi$ of the quantum affine space algebra $\AA_\qb$ 
satisfying \eqref{q-aff} there exist $u_1, \ldots, u_N \in Z(\TT_q)^{\geq 1}$
such that $u_k X_k \in \AA_\qb$ and $\psi(X_k) = (1 + u_k) X_k$, 
$\forall k \in [1,N]$.  
\eco

Artamonov \cite{Ar} studied the automorphisms of completions 
of quantum tori 
$\TT_\qb$ with respect to maximal valuations 
$\vv \colon \TT_\qb \backslash \{ 0 \} \to \Zset^N$ 
in a different direction from ours.  He considers
quantum tori for which the parameters $\{q_{kl} \mid 1 \leq k < l \leq N\}$
form a free subgroup of $\KK^*$ of rank $N(N-1)/2$ 
(such tori appear more rarely and have trivial centers) and deals with all 
automorphisms as opposed to a special subclass 
of the unipotent ones.
 
For the rest of the section
we will use the notation in the left hand side of \eqref{e-gen} for 
the generators of $\TT_\qb$,
which is more instructive in working with quantum tori.
\subsection{Supports and restrictions of unipotent automorphisms}
\label{3.2} 
For an element 
\[
u = \sum_{f \in \ZN} c_f X^f \in \wh{\TT}_{\qb, \db}
\]
denote $[u]_f:= c_f$. Define its support by
\[
\Supp(u) := \{ f \in \ZN \mid [u]_f \neq 0 \}.
\]
\bde{Sup} Given a finite unipotent automorphism $\phi$
of $\wh{\TT}_{\qb, \db}$, we will call the set 
\begin{align*}
\Supp(\phi) &= \bigcup_{k=1}^n 
\Supp ( \phi(X^{e_k}) X^{-e_k} - 1  ) 
\\
&= \{ f \in \ZN \backslash \{ 0 \}
\mid  [\phi(X^{e_k})]_{e_k +f} \neq 0 
\; \; \mbox{for some} \; \; k \in [1,N] \}  
\end{align*}
the {\em{support}} of $\phi$.
\ede

By a strict cone $C$ in $\Rset^N$ 
we will mean a set of the 
form $\Rset_{\geq 0} X$ for a finite 
subset $X = \{x_1, \ldots, x_n \}$ of $\Rset^N$ 
such that
\begin{equation}
\label{cone} 
a_1 x_1 + \ldots + a_n x_n = 0, \; 
a_1, \ldots, a_n \in \Rp 
\Rightarrow a_1 = \ldots = a_n = 0.
\end{equation}  
In other words $C$ contains no lines. 
The support of any finite unipotent automorphism $\phi$
of $\wh{\TT}_{\qb, \db}$ lies in a translated half space
\[
\Supp(\phi) \subset \{ (a_1, \ldots, a_N) \in \Rset^N 
\mid 
d_1 a_1 + \ldots + d_N a_N \geq 1 \}.
\]
Thus for every finite collection $\phi_1, \ldots, \phi_j$
of finite unipotent automorphisms of $\wh{\TT}_{\qb, \db}$, the set 
$\Rp ( \Supp(\phi_1) \cup \ldots \cup \Supp(\phi_j) )$
satisfies the condition \eqref{cone}. 

\bde{Cone} Given a finite family $\phi_1, \ldots, \phi_j$
of finite unipotent automorphism of $\wh{\TT}_{\qb,\db}$, we define 
its joint cone by 
\[
\Con(\phi_1, \ldots, \phi_j) = 
\Rp ( \Supp (\phi_1) \cup \ldots \cup \Supp(\phi_j) )
\subset \Rset^N.
\]
\ede

\ble{Con-bifin} (i) If $\phi$ is a finite unipotent automorphism
of $\wh{\TT}_{\qb,\db}$ and $g \in \Zset^N$, 
then
\[
\Supp( \phi(X^g)X^{-g}  - 1 ) \in \big( \Zp \Supp(\phi) \big) \backslash \{ 0 \}. 
\]

(ii) If $\phi$ and $\psi$ are two finite unipotent automorphisms of $\wh{\TT}_{\qb, \db}$,
then   
\[
\Supp (\phi \ci \psi) \subset 
\big(\Zp (\Supp (\phi) \cup \Supp(\psi) ) \big) \backslash \{ 0\}
\]
and thus $\Con (\phi \ci \psi) \subset \Con(\phi) + \Con(\psi)$.

(iii) For all bi-finite unipotent automorphisms 
of $\wh{\TT}_{\qb, \db}$ we have
\[
\Supp (\phi^{-1}) \subset \big(\Zp \Supp (\phi) \big) 
\backslash \{ 0\} 
\; \; 
\mbox{and} \; \; 
\Supp (\phi) \subset \big( \Zp \Supp (\phi^{-1}) \big)
\backslash \{ 0 \}.
\] 
In particular, $\Con(\phi^{-1}) = \Con(\phi)$.
\ele
\begin{proof} (i) The coefficients of $\phi(X^{-e_k})$ are 
determined from the ones of $\phi(X^{e_k})$ using the equality
\begin{equation}
\label{inv}
\phi(X^{-e_k}) = X^{- e_k} \sum_{m=0}^\infty \big(1 - \phi(X^{e_k})X^{-e_k} \big)^m,  
\end{equation}
cf. \eqref{unit2}.
The case $g=-e_k$ of part (i) follows from this. The general case is obtained 
by multiplying such expressions for $\phi(X^{\pm e_k})$, $k \in [1,N]$.

Part (ii) is a direct consequence of part (i).

(iii) The coefficients $[\phi^{-1}(X^{e_k})]_{e_k +f}$
are recursively determined for $D(f)=1, 2, \ldots$ from the 
coefficients of $\phi$ by
\begin{align*}
[\phi^{-1}(x^{e_k})]_{e_k +f} = &- \sum_{j_1, \ldots, j_N, f_1, \ldots f_N} 
q_* [\phi(X^{e_k})]_{j_1 e_1 + \ldots + (j_k+1)e_k + \ldots + j_N e_N}
\\
&\times
[\phi^{-1}(X^{j_1 e_1})]_{j_1 e_1 + f_1}
\ldots 
[\phi^{-1}(X^{(j_k +1)e_k})]_{(j_k+1) e_k + f_k}
\ldots
[\phi^{-1}(X^{j_N e_N})]_{j_N e_N + f_N}
\end{align*}
for some appropriate elements $q_* \in \KK$, 
where the sum is over $j_1, \ldots, j_N \in \Zset$,
$f_1, \ldots, f_N \in \ZN$ such that $\sum_{k=1}^N j_k e_k \in \Supp(\phi)$,
$0< D(f_1), \ldots, D(f_N)<D(f)$,
$f_1 + \ldots +f_N + \sum_{k=1}^N j_k e_k  =f$.
In the right hand side,
using \eqref{inv} for $j' \in \Zset$ and $f' \in \ZN$ such that $D(f')>0$, 
one expresses $[\phi^{-1}(X^{j' e_k})]_{j' e_k + f'}$ in terms of 
$[\phi^{-1}(X^{e_k})]_{e_k + g}$ for $g \in \ZN$, $0< D(g)\leq D(f') < D(f)$.
(Note the finiteness of the sum in the right hand side of the above formula.)
Part (iii) of the lemma easily follows from this formula.
\end{proof}

Let $C$ be a strict cone in $\Rset^N$ 
and $x \in \Rset^N \backslash \{ 0 \}$.
The ray $\Rp x$ is 
an {\em{extremal ray}} of $C$, if $x \in C$ and 
for all $x_1, x_2 \in C$, $x_1 + x_2 \in \Rp v$ implies 
$x_1, x_2 \in \Rp x$.  
For a ray $\Rp f$ in $\Rset^N$ 
and $u \in \wh{\TT}_{\qb, \db}$ denote 
\[
u|_{\Rp f} = \sum_{g \in ( \Rp f \cap \ZN) } [u]_g X^g. 
\]

Let $\phi$ be a finite unipotent automorphism of $\wh{\TT}_{\qb,\db}$ with
\[
\phi(X^{e_k}) = (1 + u_k)X^{e_k}, 
u_k \in \TT_\qb^{\geq 1}, \; \; 
k \in [1,N].
\]
Let $\Rp f$ be an extremal ray of $\Con (\phi)$. 
It is straightforward to verify that the $N$-tuple
\[
(u_1|_{\Rp f }, \ldots, u_f|_{\Rp f}) 
\]
of elements of $\TT_\qb^{\geq 1}$
satisfies the conditions of \leref{unip}. 
Therefore it defines a finite unipotent 
automorphism of $\wh{\TT}_{\qb, \db}$ which 
will be denoted by $\phi|_{\Rp f}$. 
For $f \notin \Con (\phi)$, we set 
$\phi|_{\Rp f} : = \id$. 

The proof of the following result is analogous to the proof 
of \leref{Con-bifin} and is left to the reader. 
\bpr{restr} For all completed quantum tori $\wh{\TT}_{\qb, \db}$ and 
$f \in \ZN \backslash \{0\}$ the following hold:

(i) If $\phi$ is a finite unipotent automorphism
of $\wh{\TT}_{\qb, \db}$ and $f$ is such that 
$\Rp f$ is an extremal ray of $\Con(\phi)$ 
or $f \notin \Con(\phi)$, then
\[
\phi|_{\Rp f} (X^g) = (\phi(X^g) X^{-g})|_{\Rp f} X^g, \; \; 
\forall g \in \ZN. 
\]

(ii) If $\phi$ and $\psi$ are finite unipotent automorphisms 
of $\wh{\TT}_{\qb,\db}$ and $\Rp f$ is an extremal ray of 
$\Con(\phi, \psi)$, then 
\[
(\phi \circ \psi)|_{\Rp f} = \phi|_{\Rp f} \circ \psi|_{\Rp f},
\]
cf. \leref{Con-bifin} (ii).

(iii) If $\phi$ is a bi-finite unipotent 
automorphism of $\wh{\TT}_{\qb, \db}$ and 
$\Rp f$ is an extremal ray of $\Con(\phi)$, 
then $\phi|_{\Rp f}$ is a bi-finite unipotent 
automorphism of $\wh{\TT}_{\qb, \db}$ and 
\[
\big( \phi|_{\Rp f}\big)^{-1} = (\phi^{-1})|_{\Rp f},
\]
cf. \leref{Con-bifin} (iii).
\epr
\subsection{Bi-finite unipotent automorphisms of completed quantum tori}
\label{3.3}
Our proof of \thref{m1} is based on a result for unipotent automorphisms
of completed saturated quantum tori $\wh{\TT}_{\qb, \db}$ with support not lying 
in $\Ker (\qb)$ and on \prref{restr}. The former is obtained
in \prref{non-f}. \prref{examp1} and 
\leref{const} contain two auxiliary results for the proof of \prref{non-f}.

For $f \in \ZN$ denote by $\mu_f \colon \wh{\TT}_{\qb, \db} \to \wh{\TT}_{\qb, \db}$ 
the continuous $\KK$-linear map, given by 
\[
\mu_f (X^g) = X^f X^g, \; \; \forall g \in \ZN.
\]

\bpr{examp1} Let $\wh{\TT}_{\qb, \db}$ be a completed saturated 
quantum torus of rank $N$ over an arbitrary field $\KK$ as above.
For all $f \in \ZN$, $f \notin \Ker(\qb)$ and $c \in \KK$,
\[
\phi_{f, c} (X^g) := X^g + ( 1 - q^{-1}_{f, g}) 
\sum_{m=1}^\infty c^m (X^f)^m X^g, \; \; 
g \in \ZN,
\] 
where $q_{f,g} := X^f X^g X^{-f} X^{-g} \in \KK^*$,
defines a unipotent automorphism of $\wh{\TT}_{\qb, \db}$.
\epr
\begin{proof} Only the fact that $\phi_{f,c}$ is an 
endomorphism requires a proof since the bijectivity 
of the map is immediate, see \reref{unip-aut} (i).
It is easy to verify this statement directly. 
We give another proof that explains the origin of $\phi_{f,c}$. 
The lemma is effectively a statement for  
quantum tori defined over the ring $\Zset[c, q_{kl}^{\pm1}, 1\leq k < l \leq N]$
(where $c, q_{kl}^{\pm1}$ are independent variables).
It is sufficient to prove that $\phi_{f,c}$ is an 
endomorphism for base 
fields of characteristic 0. For such, it is obvious 
that $\psi_{f,c} \colon \wh{\TT}_{\qb, \db} \to \wh{\TT}_{\qb, \db}$,
given by 
\[
\psi_{f,c} = \exp \left( \sum_{m=1}^\infty \frac{c^m}{m} \ad_{(X^f)^m} \right),
\]
is a unipotent automorphism of $\wh{\TT}_{\qb,\db}$.
It satisfies 
\[
\psi_{f, c} (X^g) = \exp \left( \sum_{m=1}^\infty
\frac{c^m}{m} \big( 1 - q^{-m}_{f,g} \big) \mu_f^m \right) X^g. 
\]
The following identity of formal power series over $\Qset[q^{\pm 1},c]$
in the variable $z$ shows that $\psi_{f, c} = \phi_{f,c}$:
\begin{align*}
& \exp \Big( \sum_{m=1}^\infty c^m(1 - q^{-m})z^m/m \Big) 
= \exp ( - \log ( 1 - c z) + \log (1 - c q^{-1} z) )
\\
& = (1 - c q^{-1} z) ( 1 - c z)^{-1} 
= 1 + (1 - q^{-1}) \sum_{m=1}^\infty c^m z^m, \; \; 
\forall j \in \Zset,
\end{align*}
where $\log(1-z) := - \sum_{m=1}^\infty z^m/m$. 
Therefore $\phi_{f,c}$ is a unipotent automorphism 
of $\wh{\TT}_{\qb,\db}$
if the base field has characteristic $0$, and by the above 
reasoning for all base fields.
\end{proof}
\ble{const} Let $\wh{\TT}_{\qb, \db}$ be a
completed saturated quantum torus over an arbitrary field $\KK$ 
and 
$f \in \ZN$ such that $f \notin \Ker(\qb)$ and $f/j \notin \ZN$, 
$\forall j \in \Zset_+$. Then for every unipotent automorphism 
$\phi$ of $\wh{\TT}_{\qb,\db}$ with $\Supp \phi \subseteq \Zp f$ 
and $m \in \Zp$, there exist $c_1, \ldots, c_m \in \KK$ 
such that  
\[
\big( \phi(X^{g}) - \phi_{m f, c_m} \ldots \phi_{f, c_1}(X^{g})\big) X^{-g}
\in \Span \{ X^{(m+1)f}, X^{(m+2) f}, \ldots \}, \; \; 
\forall g \in \ZN. 
\]
\ele
In other words, every automorphism $\phi$ satisfying the 
conditions in \leref{const} is given by  
\[
\phi = \ldots \phi_{2f,c_2} \phi_{f,c_1}
\]
for some $c_1, c_2, \ldots \in \KK$.
Note that the infinite (right-to-left) product is well 
a defined unipotent automorphism 
of $\wh{\TT}_{\qb,\db}$.
\medskip
\\
\noindent
{\em{Proof of \leref{const}}}. We prove the statement by induction on 
$m$, the case $m=0$ being trivial. Assume its validity for some $m \in \Zp$,
and define 
\[
\delta(X^g) : = 
\big[ \phi(X^{g}) - \phi_{m f, c_m} \ldots \phi_{f, c_1}(X^{g})\big]_{g+ (m+1)f}
X^{g + (m+1)f}, \; \; g \in \ZN.
\]
Then $\delta$ is a derivation of $\TT_\qb$, which must be inner 
by \cite[Corollary 2.3]{OP} since $(m+1)f \notin \Ker (\qb)$, 
recall Eq. \eqref{c1}. Therefore 
$\delta = c_{m+1} \ad_{X^{(m+1)f}}$ for some $c_{m+1} \in \KK$. 
Then 
\[
X^{-g} \big( \phi(X^{g}) - \phi_{(m+1) f, c_{m+1}} \ldots \phi_{f, c_1}(X^{g})\big)
\in \Span \{ X^{(m+2)f}, X^{(m+3) f}, \ldots \}, \; \; 
\forall g \in \ZN
\]
which completes the induction.
\qed

\bpr{non-f} Assume that $\wh{\TT}_{\qb, \db}$ is a
completed saturated quantum torus of rank $N$ over an arbitrary base field $\KK$. 
If $\phi$ is a bi-finite unipotent automorphism of $\wh{\TT}_{\qb, \db}$ 
with support in $\Rp f$ for some $f \in \ZN$ such that $f \notin \Ker (\qb)$, 
then $\phi$ is the identity automorphism. 
\epr

Nontrivial finite unipotent automorphisms 
with support in $\Rp f$ for $f \in \ZN$, $f \notin \Ker (\qb)$ 
which are not bi-finite are constructed in \cite[Proposition 4.2]{BZ}. 
\medskip
\\
\noindent
{\em{Proof of \prref{non-f}}}. We have $\Rp f \cap \ZN = \Zp f_0$ for some $f _0 \in \ZN$.
Assume that the statement of the proposition is not correct.
Then there exists $k \in [1,N]$ such that 
\[
\phi(X^{e_k}) = X^{e_k} + c_1 X^{f_0} X^{e_k} + \ldots + c_m X^{m f_0} X^{e_k}
\]
for some $m \in \Zset_+$, $c_1, \ldots, c_m \in \KK$, $c_m \neq 0$. 
Moreover,
\[
\phi^{-1}(X^{e_k}) = X^{e_k} + c'_1 X^{f_0} X^{e_k} + \ldots + c'_j X^{j f_0} X^{e_k}
\]
for some $j \in \Zp$, $c'_1, \ldots, c'_j \in \KK$, $c'_j \neq 0$.
\leref{const} implies that $\phi^{-1}(X^{f_0}) = X^{f_0}$
and thus
\[
\phi^{-1} \phi (X^{e_k}) - c_m c'_j X^{m f_0} X^{j f_0} X^{e_k} \in 
\Span \{ X^{e_k}, \ldots, X^{e_k + (m+j-1) f_0} \}.
\]
Therefore $[\phi^{-1} \phi (X^{e_k})]_{e_k + (m+j) f_0} \neq 0$
which contradicts with $\phi^{-1} \phi = \id$ because $m>0$. 
\qed
\medskip

We now proceed with the proof of \thref{m1}.
\medskip
\\
\noindent
{\em{Proof of \thref{m1}.}} Assume that $\phi$ is a bi-finite
unipotent automorphism of $\wh{\TT}_{\qb, \db}$. If $\phi = \id$, 
we are done. 
Otherwise, let $\Rp f_1, \ldots, \Rp f_n$ be the extremal 
rays of $\Con(\phi)$. By \prref{restr} (iii),
$\phi|_{\Rp f_1}, \ldots, \phi|_{\Rp f_n}$ are 
bi-finite unipotent automorphisms of $\wh{\TT}_{\qb, \db}$
that are not equal to the identity automorphism.
\prref{non-f} implies that $f_1, \ldots, f_n \in \Ker (\qb)$. 
Therefore
\[
\Con(\phi) = \Rp f_1 + \ldots + \Rp f_n \subset \Ker(\qb) \otimes_{\Zset} \Rset. 
\]
It follows from Eq. \eqref{c1} that $\Supp(\phi) \subset \Ker (\qb)$.
The definition of $\Supp(\phi)$ and Eq. \eqref{Cent} 
imply that there exists an $N$-tuple
$(u_1, u_2, \ldots, u_N)$ of elements of $Z(\TT_\qb)^{\geq 1}$
such that $\phi(X_k) = (1+u_k) X_k$, $\forall k \in [1,N]$.
This completes the proof the theorem.
\qed
\sectionnew{Unipotent automorphisms of $\UU_q^-(\g)$} 
\label{UqUnip}
\subsection{Statement of the main result}
\label{4.1}
In this section we carry out the major step of the
proof of the Andruskiewitsch--Dumas conjecture.
We define unipotent automorphisms of the algebras 
$\UU^-_q(\g)$ in a similar fashion to the case 
of completed quantum tori, see \deref{Uq-unip} for 
details. In \thref{Uq-ut} we prove a rigidity result for 
them stating that every unipotent automorphism of 
$\UU_q^-(\g)$ is equal to the identity automorphism.
The proof of this result is obtained in several 
reduction stages using the rigidity result from 
the previous section and Theorems \ref{t1a} and 
\ref{t1b}. The reductions appear in 
\S \ref{4.2}--\ref{4.3} and the proof of \thref{Uq-ut}
is given in \S \ref{4.4}.

Every strictly dominant integral
coweight $\la = \sum_{\al \in \Pi} m_\al \vpi_\al\spcheck \in \PP_{++}\spcheck$ 
gives rise to a specialization of 
the $(-\QQ_+)$-grading of the algebra $\UU_q^-(\g)$
to an $\Zp$-grading as follows:
\begin{equation}
\label{grading}
\UU^-_q(\g) = \bigoplus_{m \in \Zp} \UU^-_q(\g)^m, \; \;
\UU_q^-(\g)^m := \bigoplus \{ \UU_q^-(\g)_{-\ga} \mid 
\ga \in \QQ_+, \lcor \la, \ga \rcor = m \}.
\end{equation}
In other words, the generators $F_\al$ of $\UU^-_q(\g)$ are 
assigned degrees $m_\al= \lcor \la , \al \rcor$. (The graded components 
in \eqref{grading} depend on the choice of $\la$, but this dependence will 
not be explicitly shown for simplicity of the notation
as it was done for quantum tori.) This $\Zp$-grading of $\UU_q^-(\g)$ 
is connected, i.e., $\UU^-_q(\g)^0 = \KK$.
For $m \in \Zp$, denote
\[
\UU_q^-(\g)^{\geq m} = \bigoplus_{j \geq m} \UU^-_q(\g)^j.
\]
For a graded subalgebra $R$ of $\UU_q^-(\g)$, set 
$R^{\geq m} = R \cap \UU_q(\g)^{\geq m}$. 
\bde{Uq-unip} Given a strictly dominant integral coweight $\la$,
we call an automorphism $\Phi$ of $\UU^-_q(\g)$ 
$\la$-{\em{unipotent}} if 
\[
\Phi(F_\al) - F_\al \in \UU^-_q(\g)^{\geq \lcor \la, \al \rcor + 1} \; \;
\mbox{for all} \; \; \al \in \Pi.
\]
\ede
Every $\la$-unipotent automorphism $\Phi$ of $\UU^-_q(\g)$ satisfies
\begin{equation}
\label{u5}
\Phi(u) - u \in \UU^-_q(\g)^{\geq m+1} \; \;
\mbox{for all} \; \; u \in \UU^-_q(\g)^{m}, m \in \Zp.
\end{equation}
The set of $\la$-unipotent automorphisms of $\UU_q(\g)$ 
is a subgroup of $\Aut( \UU_q^-(\g))$ for all $\la \in \PP_{++}\spcheck$.

The following theorem is the major step in our proof of the Andruskiewitsch--Dumas
conjecture.

\bth{Uq-ut} Let $\g$ be a simple Lie algebra of rank $r>1$, $\KK$ 
an arbitrary base field, and $q$ a deformation parameter that is 
not a root of unity. For every 
strictly dominant integral coweight $\la$,
the only $\la$-unipotent automorphism 
of $\UU^-_q(\g)$ is the identity automorphism.
\eth

An analogous statement holds for the algebras $\UU_q^+(\g)$ 
because of the isomorphism $\om \colon \UU_q^\pm(\g) \to \UU_q^\mp(\g)$,
see \eqref{om}. \thref{Uq-ut} trivially holds in the case when $\rank(\g)=1$ 
since $\UU_q^-(\sl_2) = \KK[F_\al]$.
\subsection{Relations to automorphisms of completed quantum tori}
\label{4.2}
Let 
\begin{equation}
\label{ib}
\ib = (\al_1, \ldots, \al_N)
\end{equation}
be a reduced word for the longest element $w_0$ of $W$. 
Recall from \S \ref{2.3} that we denote by  
$(\ol{F}_{\ib, 1}, \ldots, \ol{F}_{\ib, N})$ the 
final $N$-tuple of elements of $\Fract(\UU^-_q(\g))$ 
from the Cauchon deleting derivation procedure 
for the iterated Ore extension presentation \eqref{Umin-iter} 
of $\UU^-_q(\g)$. Denote by $\TT(\ib)$ the subalgebra 
of $\Fract(\UU^-_q(\g))$ generated by them and their inverses
\[
\TT(\ib) = \KK \lcor \ol{F}_{\ib, 1}^{\, \pm 1}, \ldots, \ol{F}_{\ib, N}^{\, \pm 1} 
\rcor \subset \Fract ( \UU^-_q(\g)),
\]
see \eqref{C-tor}. 
Let $\qb \in M_N(\KK^*)$ be the multiplicatively skew-symmetric 
matrix such that  
\begin{equation}
\label{q-matr}
q_{kl} = q^{ \lcor \be_k, \be_l \rcor} \; \; 
\mbox{for} \; \; 1 \leq k < l \leq N.
\end{equation}
By \eqref{Cau} and \eqref{coeff} we have the isomorphism of quantum tori
\[
\TT_\qb \cong \TT(\ib), \; \; 
X_l \mt \ol{F}_{\ib, l}, \; \; 
l \in [1,N]
\]
in the notation of \eqref{q-tor}. We change the generating 
set $(\ol{F}_{\ib,1}, \ldots, \ol{F}_{\ib,N})$ 
using \thref{1b}. Recall the definition \eqref{succ} of the successor 
function $s \colon [1,N] \to [1,N] \cup \{ \infty \}$ associated 
to the reduced word $\ib$. For $l \in [1,N]$ denote
\[
O(l) = \max \{ m \in \Zp \mid s^m(j) \neq \infty \},
\]
where as usual $s^0 = \id$. \thref{1b} implies 
\begin{equation}
\label{main1b-eq}
\De_{\ib, l} = (q_{\al_l}^{-1} - q_{\al_l})^{O(l)} 
\ol{F}_{\ib,s^{O(l)}(l)} \ldots \ol{F}_{\ib,l}
\end{equation}
and
\[
\TT(\ib) = \KK \lcor \De_{\ib, 1}^{\, \pm 1}, \ldots, \De_{\ib, N}^{\, \pm 1} 
\rcor \subset \Fract ( \UU^-_q(\g)).
\]
Moreover, from \eqref{main1b-eq} we have
\[
\De_{\ib, k} \De_{\ib, l} = 
q^{n'_{kl}} \De_{\ib, l} \De_{\ib, k},
\]
where
\begin{equation}
\label{nns}
n'_{kl} = \sum_{j=1}^{O(k)} \sum_{m=1}^{O(l)} 
\sign( s^m(l) - s^j(k)) \lcor \be_{s^j(k)}, \be_{s^m(l)} \rcor
\; \; \mbox{for} \; \; k, l \in [1, N].
\end{equation}
Denote by $\qb' \in M_N(\KK)$ the multiplicatively skew-symmetric matrix
whose entries are given by $q'_{kl} = q^{n'_{kl}}$.
We have the isomorphism
\begin{equation}
\label{imp}
\TT_{\qb'} \cong \TT(\ib), \; \; 
X_l \mt \De_{\ib, l}, \; \; 
l \in [1,N].
\end{equation}
The $\QQ$-grading of $\UU^-_q(\g)$ gives rise to 
a $\QQ$-grading of the quantum torus $\TT(\ib)$ by assigning
\[
\deg \ol{F}_{\ib, l} := \deg F_{\be_l} = - \be_l.
\]
Given a strictly dominant integral coweight $\la \in \PP_{++}\spcheck$, 
this grading can be specialized to a $\Zset$-grading 
\begin{equation}
\label{Tgrad}
\TT(\ib) = \bigoplus_{m \in \Zset} \TT(\ib)^m
\; \; \mbox{by setting} 
\; \; \deg \ol{F}_{\ib,l} := \lcor \la, \be_l \rcor.
\end{equation}
This grading is compatible 
with the $\Zp$-grading \eqref{grading} of $\UU^-_q(\g)$
\begin{equation}
\label{incl-m}
\UU^-_q(\g)^m \subset \TT(\ib)^m, \; \; 
\forall m \in \Zp,
\end{equation}
recall \eqref{Cau}.
Because of \eqref{main1b-eq},  
the $\Zset$-grading of $\TT(\ib)$ can be also 
defined in terms of the generators $\De_{\ib,l}$
by 
\[
\deg \De_{\ib, l} = \sum_{j=0}^{O(l)} \lcor \la, \be_{s^j(l)} \rcor.
\]
We associate to $\la \in \PP_{++}\spcheck$ the following degree vector
\begin{equation}
\label{deg-vect}
\db = (d_1, \ldots, d_N) \in \Zset_+^N,
\; \; 
d_l := \sum_{j=0}^{O(l)} \lcor \la, \be_{s^j(l)} \rcor, \; \; 
l \in [1,N]
\end{equation}
for the quantum torus $\TT(\ib) \cong \TT_{\qb'}$, recall \S \ref{3.1}.
The grading \eqref{Tgrad} is precisely the grading associated 
to $\db$ as defined in \S \ref{3.1}.
Denote by $\wh{\TT}(\ib, \db)$ the corresponding completion of $\TT(\ib)$, 
which is isomorphic to $\wh{\TT}_{\qb',\db}$. 

{\em{To every $\la$-unipotent automorphism of $\UU^-_q(\g)$ we 
will associate a bi-finite unipotent automorphism of 
the completed quantum torus $\wh{\TT}(\ib, \db)$ as follows}},
recall Definitions \ref{dunip} and \ref{dUq-unip}.
Let $\Phi$ be a unipotent automorphism of $\UU^-_q(\g)$.
By Eqs. \eqref{u5} and \eqref{incl-m}
\[
\Phi(\De_{\ib, l}) = \De_{\ib, l} + u'_l \; \; 
\mbox{for some} \; \; 
u'_l \in \UU^-_q(\g)^{\geq d_l + 1}
\subseteq \TT(\ib)^{\geq d_l + 1},
l \in [1,N].
\]
Therefore 
\[
\Phi(\De_{\ib, l}) = (1 + u_l) \De_{\ib, l} \; \; 
\mbox{for} \; \; 
u_l := u'_l \De_{\ib, l}^{-1} \in \TT(\ib)^{\geq 1}, 
l \in [1,N].
\]
Since $\Phi$ is an automorphism of $\UU_q^-(\g)$, we obtain from \eqref{imp} 
that the 
$N$-tuple $(u_1, \ldots, u_N)$ satisfies the condition
in \leref{unip}. Therefore there exists a unique finite 
unipotent automorphism $\phi$ of $\wh{\TT}(\ib,\db)$ given by 
\[
\phi(\De_{\ib, l}) := \Phi(\De_{\ib, l}) = (1 + u_l) \De_{\ib, l},
\; \; \forall l \in [1,N]. 
\]  
We have $\UU^-_q(\g) \subset \TT(\ib) \subset \wh{\TT}(\ib, \db)$ and 
\begin{equation}
\label{restr}
\Phi|_{\UU_q^-(\g)} = \phi.
\end{equation}
Denote by $\psi$ the finite unipotent automorphism of $\wh{\TT}(\ib,\db)$
associated to the inverse unipotent automorphism $\Phi^{-1}$ of 
$\UU^-_q(\g)$ by the above construction.
It follows from \eqref{restr} that  
$(\phi \circ \psi)|_{\UU_q^-(\g)} = (\psi \circ \phi)|_{\UU_q^-(\g)} = \id$. 
Because $\De_{\ib, l} \in \UU_q^-(\g)$, $\forall l \in [1,N]$
and those elements generate $\TT(\ib)$,  
we have $\phi^{-1} = \psi$. Thus, this construction associates a 
bi-finite unipotent automorphism $\phi$ of $\wh{\TT}(\ib, \db)$
to each $\la$-unipotent automorphism $\Phi$ of $\UU_q^-(\g)$.

\bre{Ca} We note that one could not achieve the same result using 
the Cauchon elements $\ol{F}_{\ib,1}, \ldots, \ol{F}_{\ib, N} \in \Fract(\UU_q^-(\g))$.
The key point is that they and their inverses do generate a 
quantum torus inside $\Fract(\UU_q^-(\g))$, but the Cauchon 
elements do not belong to $\UU_q^-(\g)$ in general (i.e. 
the first inclusion in \eqref{me} is not satisfied for the 
quantum affine space algebra generated by the Cauchon elements).
Because of this 
they cannot be used for the passage to bi-finite unipotent automorphisms 
of $\wh{\TT}(\ib, \db)$. What we 
did was to change the Cauchon 
generators $\ol{F}_{\ib,1}^{\pm 1}, \ldots, \ol{F}_{\ib, N}^{\pm 1}$
of $\TT(\ib)$ to a new generating set 
$\De_{\ib, 1}^{\pm 1}, \ldots, \De_{\ib, N}^{\pm 1}$ of $\TT(\ib)$ 
in such a way that $\De_{\ib, 1}, \ldots, \De_{\ib, N} \in \UU^-_q(\g)$.
Because of these two properties, this new set of generators can be used to establish the connection to bi-finite 
unipotent automorphisms of completed quantum tori as indicated.
\ere

We use the above relationship and \thref{m1} to obtain our first reduction 
for the proof of \thref{Uq-ut}. 

Recall that the longest element $w_0$ of the Weyl group $W$
gives rise to the involution
\begin{equation}
\label{w0}
\al \mt \wt{\al} := - w_0(\al)
\end{equation}
of $\Pi$, which is an element of $\Aut(\Ga)$. 

\bpr{ste1} Let $\g$, $\KK$, $q \in \KK^*$,
and $\la \in \PP_{++}\spcheck$ be as in \thref{Uq-ut}.
Assume that $\al \in \Pi$ and $\ib=(\al_1, \ldots, \al_N)$ 
is a reduced word for $w_0$ such that 
\begin{equation}
\label{al-N}
\al_N = \wt{\al}.
\end{equation}
If $\Phi$ is a $\la$-unipotent automorphism of $\UU^-_q(\g)$, 
then there exists
\[
z_\al \in Z(\TT(\ib))^{\geq 1}
\]  
such that 
\[
\Phi(F_\al) = (1 + z_\al) F_\al.
\]
\epr
\begin{proof} Since $q \in \KK^*$ is not a root of unity, the subgroup of 
$\KK^*$ generated by all entries of the matrix $\qb'$ is 
torsion-free. Therefore the quantum tori $\TT(\ib)$ and $\TT_{\qb'}$
are saturated, recall \eqref{imp}. We apply \thref{m1} to 
the bi-finite unipotent automorphism $\phi$ of $\wh{\TT}(\ib,\db)$
associated to $\Phi$. It implies that
\[
\Phi(\De_{\ib,N}) \De_{\ib,N}^{-1} = 
\phi(\De_{\ib,N}) \De_{\ib,N}^{-1} = 1 + z_\al
\]
for some $z_\al \in Z(\TT(\ib))^{\geq 1}$.
\thref{1b} and Eq. \eqref{new-x} imply
\[
(q_{\al_N}^{-1} - q_{\al_N})^{-1} \De_{\ib,N} = \ol{F}_{\ib, N} = F_{\be_N}.
\] 
Since $\be_N = s_1 \ldots s_{N-1}(\al_N) = - w_0 s_{\al_N} (\al_N) = \wt{\al}_N 
= \al$, it follows from Eq. \eqref{eqa} that 
\begin{equation}
\label{beal}
F_{\be_N} = F_\al.
\end{equation}
Therefore $\Phi(F_\al) = (1+ z_\al) F_\al$. 
\end{proof}
\bre{short} One can rephrase the arguments from the first part of this subsection 
(the proof of Eq. \eqref{imp}) in a shorter form, which only uses \cite[Proposition 3.3]{GY} 
and \cite[Theorem 10.1]{BZ}, and does not use \cite{Ca}. 
This proof has the disadvantage that the passage from unipotent automorphisms 
of $\UU_q^-(\g)$ to bi-finite unipotent automorphisms of completed quantum tori 
misses the conceptual point described in \reref{Ca}.
For the sake of completeness we sketch this shorter route.
It follows from Theorem 10.1 in \cite{BZ} that 
\begin{equation}
\label{Decomm}
\De_{\ib, k} \De_{\ib, l} = 
q^{n'_{kl}} \De_{\ib, l} \De_{\ib, k}, \quad \forall k, l \in [1,N]
\end{equation}
for some $n'_{kl} \in \Zset$, see \cite[\S 3.3]{GY} for details. 
By \cite[Proposition 3.3]{GY} for all $l \in [1,N]$,
\begin{align*}
&\De_{\ib, l} = 
(q_{\al_l}^{-1} - q_{\al_l}) \De_{\ib, s(l)} F_{\be_l}
\mod \UU^-_q(\g)_{\ib, [l+1, N]},
\quad \quad \mbox{if} \; \; s(l) \neq \infty, 
\\
&\De_{\ib, l} = 
(q_{\al_l}^{-1} - q_{\al_l}) F_{\be_l}
\mod \UU^-_q(\g)_{\ib, [l+1, N]}, 
\quad \quad \quad \quad \; \; \mbox{if} \; \; s(l)= \infty,
\end{align*}
where $\UU^-_q(\g)_{\ib, [l+1, N]}$ denotes the unital subalgebra of $\UU_q^-(\g)$
generated by $F_{\be_{l+1}}, \ldots, F_{\be_N}$. (Recall \S \ref{2.4} 
for the definition of $s \colon [1,N] \to [1,N] \sqcup \{\infty\}$.)
The above three facts 
and the Levendorskii--Soibelman straightening law \eqref{LS} imply 
that the $\KK$-subalgebra of $\Fract(\UU^-_q(\g))$ generated by  
$\De_{\ib, 1}^{\pm 1}, \ldots, \De_{\ib, N}^{\pm 1}$ 
is isomorphic to a quantum torus and that
the integers $n'_{kl}$ are given by \eqref{nns}.
The second fact follows by comparing leading terms 
in \eqref{Decomm} with respect to PBW bases and using \eqref{LS}.
This establishes the isomorphism from \eqref{imp} and the embedding
$\UU^-_q(\g) \subset \TT(\ib)$, 
and we can proceed
with the other arguments as described above. In the case when the 
base field $\KK$ has characteristic 0 and $q \in \KK$ is transcendental over $\Qset$ 
the isomorphism from \eqref{imp} and the embedding $\UU^-_q(\g) \subset \TT(\ib)$
also follow from the work of Gei\ss, Leclerc and Schr\"oer \cite{GLS}.
\ere
\subsection{Second reduction step for \thref{Uq-ut}}
\label{4.3}
Consider the involution \eqref{w0} of $\Pi$. Denote its 
fixed point set by $\Pi^0 = \{ \al \in \Pi \mid - w_0(\al) = \al \}$. 
Choose a set of base points $\Pi^+$ of its 2-element 
orbits. Let
$\Pi^- = - w_0 (\Pi^+)$. Then we have the decomposition 
\[
\Pi = \Pi^0 \sqcup \Pi^+ \sqcup \Pi^-.
\] 
The kernel 
$\Ker (1+ w_0) := \{ \la \in \PP \mid (1+ w_0) \la = 0 \}$
is given by
\begin{equation}
\label{ker}
\Ker (1 + w_0)=
\Zset \{ \vpi_\al \mid \al \in \Pi^0 \} \oplus 
\Zset \{ \vpi_{\al} + \vpi_{\wt{\al}} \mid \al \in \Pi^+ \}.  
\end{equation}
It follows from the second statement in \thref{1a} (ii) that the subalgebra 
\[
\CC_q^-(\g) = \KK \lcor b^{\vpi_{\al'}}_1, b_1^{\vpi_{\al}} b_1^{\vpi_{\wt{\al}}} \mid
\al' \in \Pi^0, \al \in \Pi^+ \rcor 
\subset \NN_q^-(\g) 
\]
is a polynomial algebra over $\KK$ in the generators 
\begin{equation}
\label{set2}
\{b^{\vpi_{\al'}}_1 \mid \al' \in \Pi^0 \} \sqcup 
\{ b_1^{\vpi_{\al}} b_1^{\vpi_{\wt{\al}}} \mid \al \in \Pi^+ \}.
\end{equation}
We show in \leref{43a} that $\CC_q^-(\g) = Z(\UU_q^-(\g))$. 
The following is the second reduction step in our proof of \thref{Uq-ut}:

\bpr{2re} Let $\g$, $\KK$, $q \in \KK^*$, and $\la \in \PP_{++}\spcheck$
be as in \thref{Uq-ut}. 
If $\Phi$ is a $\la$-unipotent automorphism of $\UU_q^-(\g)$, 
then there exist elements
\[
z_\al \in \CC_q^-(\g)^{\geq 1} \; \; 
\mbox{for} \; \; \al \in \Pi
\]
such that 
\[
\Phi(F_\al) = (1 + z_\al) F_\al, \; \; 
\forall \al \in \Pi.
\]
\epr

Before we proceed with the proof of \prref{2re} we obtain
several auxiliary results. Denote the multiplicative 
subset 
\[
\Om(\g) = \KK^* \{ b^\la_1 \mid \la \in \PP_+ \}
\]
consisting of normal elements of $\UU^-_q(\g)$,
recall \eqref{norm}.
Consider the localizations
\[
\NN^-_q(\g)^\sharp := 
\NN^-_q(\g) [ \Om(\g)^{-1}] \subset 
\UU^q_-[ \Om(\g)^{-1}]. 
\]
The second statement in \thref{1a} (ii) implies that $\NN_q^-(\g)^\sharp$ 
is a Laurent polynomial algebra over $\KK$ in 
the generators $\{ b^{\vpi_\al}_1 \mid \al \in \Pi \}$.
Denote by $\CC_q^-(\g)^\sharp$ the localization 
of $\CC_q^-(\g)$ by the multiplicative subset generated 
by the elements \eqref{set2}. 
This localization is a Laurent polynomial 
algebra over $\KK$ in the generators \eqref{set2}.

\ble{43a} In the above setting
\begin{equation}
\label{ZZ}
Z(\UU^-_q(\g) [\Om(\g)^{-1}]) = \CC^-_q(\g)^\sharp 
\end{equation}
and 
\begin{equation}
\label{Z}
Z(\UU^-_q(\g)) = \CC^-_q(\g).
\end{equation}
\ele

\leref{43a} was obtained by Caldero \cite{Cal1} in the case when $\KK = \Cset(q)$.

\begin{proof}
Each $\mu \in \PP$ can be uniquely represented as $\mu = \la_+ - \la_-$, 
where $\la_\pm \in \PP_+$ and $\Supp \la_+ \cap \Supp \la_- = \emptyset$.
Denote 
\begin{equation}
\label{b-mu-elem}
b^\mu_1 := (b^{\la_+}_1)^{-1} b^{\la_-}_1 \in \NN_q^-(\g)^\sharp.
\end{equation}
By \eqref{norm}, for all $\mu \in \PP$
\[
b_1^\mu u
= 
q^{\lcor (1+ w_0)\mu, \ga \rcor }
u b_1^\mu, \; \; 
\forall \, 
u \in \UU_q^-(\g)_{-\ga},
\ga \in \QQ_+.
\]
The first statement in \thref{1a} (ii) implies that 
$Z(\UU^-_q(\g) [\Om(\g)^{-1}])$ is equal to the span of 
those $b^\mu_1$ such that $\mu \in \Ker(1+w_0)$. Eq. \eqref{ZZ} 
now follows from Eq. \eqref{ker}. By  
\thref{1a} (iii), $\UU^-_q(\g)$ is a 
free left $\NN^-_q(\g)$-module in which $\NN^-_q(\g)$ 
is a direct summand. This implies 
$\CC^-_q(\g)^\sharp \cap \UU_q^-(\g) = \CC^-_q(\g)$
and the validity of \eqref{Z}. 
\end{proof}

Given $\al \in \Pi$ and a reduced word $\ib$ for $w_0$ as in \eqref{ib},
denote
\[
l(\al,\ib) = \min \{ l \in [1,N] \mid \al_l = \al \}.
\]
From the definition \cite[Eq. 8.6 (2)]{Ja} of the braid group 
action on $V(\vpi_\al)$ it follows at once that 
$T_{w_0(\ib)_{\leq (l(\al, \ib)-1)}^{-1}} v_{\vpi_\al} = v_{\vpi_\al}$
which implies
\[
\De_{\ib, l(\al, \ib)} = b^{\vpi_\al}_1, \forall \al \in \Pi,
\]
recall \S \ref{2.3} for notation.
Therefore for all reduced words $\ib$ for $w_0$, we have 
\[
\CC_q^-(\g)^\sharp \subset \NN_-^q(\g)^\sharp \subset 
\UU^-_q(\g) [\Om(\g)^{-1}]
\subset \TT(\ib) \subset \Fract(\UU_q(\g)).
\]
Let $z$ be a Laurent monomial in $\De_{\ib, 1}, \ldots, \De_{\ib, N}$ which belongs to the center 
of $\TT(\ib)$. Define the ideal 
\[
I:= \{ r \in \UU^-_q(\g)[\Om(\g)^{-1}] \mid z r \in \UU^-_q(\g)[\Om(\g)^{-1}] \}
\]
of $\UU^-_q(\g)[\Om(\g)^{-1}]$. It is nonzero since $\De_{\ib, 1}, \ldots, \De_{\ib, N} \in \UU^-_q(\g)$. 
Furthermore, $I$ is $\Tset^r$-stable since $z$ is 
an eigenvector of $\Tset^r$ and the base field $\KK$ is infinite. 
By \cite[Theorem 3.1 (b)]{Y6}, $\UU^-_q(\g) [\Om(\g)^{-1}]$ is a 
$\Tset^r$-simple algebra, so $I= \UU^-_q(\g) [\Om(\g)^{-1}]$ 
and $z \in Z(\UU^-_q(\g) [\Om(\g)^{-1}])$. 
It follows from \eqref{Cent} that
\[
Z(\UU^-_q(\g) [\Om(\g)^{-1}]) = 
Z(\TT(\ib)).
\]
\leref{43a} implies 
\begin{equation}
\label{ZZZ}
Z(\TT(\ib)) = \CC_q^-(\g)^\sharp \; \; 
\mbox{for all reduced words $\ib$ for $w_0$}.
\end{equation}
This fact was earlier obtained by Bell and Launois in \cite[Proposition 3.3]{BL} 
using a result of De Concini and Procesi \cite[Lemma 10.4 (b)]{DP}
describing $\Ker(\qb)$ for the matrix \eqref{q-matr}, recall \eqref{Cent}.
The advantage of the above proof is that it trivially extends to the multiparameter case, 
see Section \ref{multi}.

\ble{dir-sum} In the above setting, for all $\al \in \Pi$, the module
$\NN^-_q(\g) F_\al$
is a direct summand of $\UU^-_q(\g)$ considered as a 
left $\NN^-_q(\g)$-module. In particular, if 
$u \in \NN^-_q(\g)^\sharp$ and $u F_\al \in \UU^-_q(\g)$, 
then $u \in \NN^-_q(\g)$. 
\ele
\begin{proof} Let $\ib$ be a reduced word for $w_0$ as in 
\eqref{ib} such that $\al_N= \wt{\al}$. This implies 
$F_{\be_N} = F_\al$, see \eqref{beal}. The first part 
of the lemma will follow from \thref{1a} (iii) once we show
\[
(0, \ldots, 0, 1) \in \HH_{\ib},
\]
recall \eqref{Hset} for the definition of the 
set $\HH_\ib$. This is in turn equivalent to saying
that there exists $l < N$ such that $\al_l = \al_N$. Assume 
the opposite. Then $w_0 s_{\al_N}$ belongs to the 
parabolic subgroup $W_{\al_N}$ of $W$ generated by
all simple reflections except $s_{\al_N}$. This is 
a contradiction since it implies that the longest 
element of $W_{\al_N}$ has length $N-1$, which is 
impossible if $r = \rank(\g)>1$.
It proves the first statement of the lemma. The second 
statement is a direct consequence of the first one.
\end{proof} 
\noindent
{\em{Proof of \prref{2re}}}. Let $\Phi$ be a $\la$-unipotent 
automorphism of $\UU^-_q(\g)$. By \prref{ste1} and 
Eq. \eqref{ZZZ} there exist elements
\[
z_\al \in (\CC^-_q(\g)^\sharp)^{\geq 1} \; \; 
\mbox{for} \; \; \al \in \Pi
\]
such that $\Phi(F_\al) = (1 + z_\al) F_\al$, $\forall \al \in \Pi$. 
\leref{dir-sum} implies that 
$z_\al \in \CC^-_q(\g)^\sharp \cap \NN_q^-(\g) = \CC^-_q(\g)$
because $\NN_q^-(\g)$ is a polynomial algebra.
Hence 
\[
z_\al \in \CC^-_q(\g)^{\geq 1}, \; \; \forall \al \in \Pi,
\]
which completes the proof of the proposition.
\qed
\subsection{Proof of \thref{Uq-ut}}
\label{4.4}
In this subsection we complete the proof of the triviality 
of all unipotent automorphisms of the algebras $\UU_q^-(\g)$.
\medskip
\\
\noindent
{\em{Proof of \thref{Uq-ut}}}.
Let $\la$ be a strictly dominant integral coweight 
and $\Phi$ a $\la$-unipotent automorphism of $\UU_q^-(\g)$. 
\prref{2re} implies that there exist elements
\[
z_\al \in \CC^-_q(\g)^{\geq 1} \; \; 
\mbox{for} \; \; \al \in \Pi
\]
such that $\Phi(F_\al) = (1 + z_\al) F_\al$, $\forall \al \in \Pi$.
For $\ga = \sum_\al m_\al \al \in \QQ_+$ denote 
\begin{equation}
\label{zga}
z_\ga := \prod_{\al \in \Pi} (1 + z_\al)^{m_\al} -1 \in \CC^-_q(\g)^{\geq 1}. 
\end{equation}
Since $\UU_q^-(\g)$ is generated by the set $\{ F_\al \mid \al \in \Pi\}$, 
we have 
\[
\Phi(u) = (1 + z_\ga) u, \; \; 
\forall \;
x \in \UU_q^-(\g)_{- \ga},
\ga \in \QQ_+.
\]
In particular, 
\[
\Phi(b^{\vpi_\al}_1) = (1 + z_{(1-w_0)\vpi_\al}) b^{\vpi_\al}_1, 
\; \; \forall \al \in \Pi,
\]
recall \eqref{bel}. Therefore 
\begin{equation}
\label{subset}
\Phi (\UU_q^-(\g) b^{\vpi_\al}_1) \subseteq \UU_q^-(\g) b^{\vpi_\al}_1, 
\; \; \forall \al \in \Pi.
\end{equation}
By \thref{1a} (iii), $\UU_q^-(\g) b^{\vpi_\al}_1$ are height one prime 
ideals of $\UU_q^-(\g)$. Hence we have equalities in 
\eqref{subset} because $\Phi$ is an automorphism.   
This implies that there exist elements
$x_\al \in \UU^-_q(\g)$ for $\al \in \Pi$ such that
\[
b^{\vpi_\al}_1 = x_\al (1 + z_{(1-w_0)\vpi_\al})
b^{\vpi_\al}_1, \; \; \forall \al \in \Pi.
\]
Furthermore, the algebras $\UU^-_q(\g)$ are domains whose groups 
of units are reduced to scalars because they are iterated Ore extensions.
Thus $u_\al (1+ z_{(1-w_0)\vpi_\al}) = 1$ and $1 + z_{(1-w_0)\vpi_\al} \in \KK^*$, 
$\forall \al \in \Pi$. It follows from Eq. \eqref{zga} that 
$z_{(1-w_0)\vpi_\al} = 0$, $\forall \al \in \PP_+$. 
From Eq. \eqref{zga} we obtain 
\[
z_{(1-w_0)\mu} = 0, \; \; \mbox{i.e.}, \; \; 
\Phi|_{\UU^-_q(\g)_{-(1-w_0)\mu}} = \id, \; \; \forall \mu \in \PP_+.
\]
Choose $\mu$ to be equal to the highest root
$\sum_{\al \in \Pi} m_\al \al \in \PP_+$ of $\g$. 
Taking into account that $(1-w_0)\mu = 2 \mu$ 
leads to 
\[
\prod_{\al \in \Pi} (1 +z_\al)^{2 m_\al} = 1.
\]
If the highest term of $z_\al$ with respect to the $\Zp$-grading 
\eqref{grading} of $\UU_q^-(\g)$ is in degree $n_\al$, then the 
left hand side has a nontrivial component in degree 
$\sum_{\al \in \Pi} 2 m_\al n_\al$. The right hand side lies in 
degree 0. Therefore $n_\al = 0$ for all $\al \in \Pi$ because 
$m_\al >0$, $\forall \al \in \Pi$. Hence we obtain $z_\al = 0$, 
$\forall \al \in \Pi$ and $\Phi= \id$.
\qed
\sectionnew{Proof of the Andruskiewitsch--Dumas conjecture} 
\label{sAD}
\subsection{Statement of the main result}
\label{4a.1}
Here we complete the proof of the 
Andruskiewitsch--Dumas conjecture. This result 
is stated in \thref{ADc}. Its proof relies on \thref{Uq-ut}
and on a classification result for a type of automorphisms of 
$\UU_q^-(\g)$ which we call {\em{linear}}, see \deref{lin}.

Recall that $\Ga$ denotes the Dynkin diagram of $\g$. The 
automorphism group of the directed graph $\Ga$
is denoted by $\Aut(\Ga)$. 
One has the embeddings
\[
\Upsilon^\pm \colon \Tset^r \rtimes \Aut (\Ga) \hra \Aut( \UU_q^\pm(\g)).
\]
To the pair $(t, \theta) \in \Tset^r \rtimes \Aut (\Ga)$, 
where $t = (t_\al)_{\al \in \Pi}$, one associates 
the automorphisms $\Upsilon_{(t, \theta)}^\pm \in \Aut(\UU_q^\pm(\g))$
given by 
\begin{align}
\label{Upsilon1}
\Upsilon_{(t, \theta)}^+ (E_\al) &= t \cdot E_{\theta(\al)} = t_{\theta(\al)} 
E_{\theta(\al)},
\\
\label{Upsilon2}
\Upsilon_{(t, \theta)}^- (F_\al) &= t \cdot F_{\theta(\al)} = t_{\theta(\al)}^{-1} 
F_{\theta(\al)} 
\end{align}
for the $\Tset^r$-action \eqref{torus-act}. 
The following theorem proves the Andruskiewitsch--Dumas conjecture:

\bth{ADc} For all simple Lie algebras $\g$ of rank $r>1$, 
base fields $\KK$, and deformation parameters $q \in \KK^*$ that are 
not roots of unity, the maps  
\[
\Upsilon^\pm \colon \Tset^r \rtimes \Aut (\Ga) \to \Aut (\UU_q^\pm(\g))
\]
are group isomorphisms.
\eth
The key point of the theorem is the surjectivity of 
the maps $\Upsilon^\pm$. Because of the isomorphism 
$\om$ from \eqref{om} the plus and minus cases 
are equivalent. For $\rank(\g) =1$, one has
$\UU_q^-(\g) = \KK[F_\al]$ and 
$\Aut (\UU_q^-(\g)) \cong \KK \rtimes \KK^*$.
The theorem is proved in 
\S \ref{4a.3}. First, we consider the following special 
type of automorphisms of $\UU_q^-(\g)$.

\bde{lin} We call an automorphism $\Phi$ of $\UU_q^-(\g)$ 
{\em{linear}} if 
\[
\Phi(F_\al) \in \Span \{ F_{\al'} \mid \al' \in \Pi \}, \; \; 
\forall \al \in \Pi.
\]
\ede
The set of all linear automorphisms of $\UU^-_q(\g)$ 
is a subgroup of $\Aut(\UU^-_q(\g))$.
In the next subsection we prove the following classification 
of linear automorphisms of $\UU_q^-(\g)$.
\bpr{line} Let $\g$ be a simple Lie algebra, $\KK$ and arbitrary base
field, and $q \in \KK^*$ not a root of unity. All
linear automorphism of $\UU^-_q(\g)$ are of the 
form $\Upsilon^-_{(t,\theta)}$ for some $t \in \Tset^r$ and
$\theta \in \Aut (\Ga)$. 
\epr
\subsection{Linear automorphisms of $\UU^-_q(\g)$}
\label{4a.2}
Before we proceed with the proof of \prref{line}, we establish 
an auxiliary lemma. For a linear automorphism $\Phi$ 
of $\UU_q^-(\g)$ given by 
$\Phi(F_\al) = \sum_{\al' \in \Pi} c_{\al\al'} F_{\al'}$, 
$c_{\al \al'} \in \KK$,
denote 
\begin{equation}
\label{chi}
\chi(\Phi, \al) = \{ \al' \in \Pi \mid c_{\al \al'} \neq 0 \}.
\end{equation}
\ble{lin-le} Assume that, in the setting of \prref{line}, $\Phi$
is a linear automorphism of $\UU_q^-(\g)$. Then the following hold:

(i) For all $\al, \al' \in \Pi$ such that 
$a_{\al \al'} = -1$ (recall \eqref{a})
we have 
\[
\chi(\Phi, \al) \cap \chi(\Phi, \al') = \emptyset \; \;
\mbox{and} \; \; 
|\chi(\Phi, \al)| = 1. 
\]

(ii) If there exist an element $\theta$ of the symmetric group $S_{\Pi}$ and  
scalars $t'_\al \in \KK^*$ for $\al \in \Pi$ such that
\begin{equation}
\label{per}
\Phi(F_\al) = t'_\al F_{\theta(\al)}, \; \; \forall \al \in \Pi,
\end{equation}
then $\theta \in \Aut(\Ga)$ and $\Phi = \Upsilon^-_{(t,\theta)}$,
where $t= (t_\al)_{\al \in \Pi} \in \Tset^r$ is given by
$t_\al = t_{\theta^{-1}(\al)}^{-1}$. 
\ele
Denote by $\FF^-(\g)$ the free $\KK$-algebra 
in the generators $\{F_\al \mid \al \in \Pi \}$. 
It is $(- \QQ_+)$-graded by assigning weight $-\al$ to $F_\al$.
For $\al \neq \al' \in \Pi$ denote by $R_{\al, \al'}$ the 
the expression in the left hand side of Eq. \eqref{S1} 
considered as an element of $\FF^-(\g)$.
Let $\II^-_q(\g)$ be the (graded) two sided ideal of $\FF^-(\g)$ 
generated by all such elements $R_{\al, \al'}$.
Denote by $\II^-_q(\g)_{- \ga}$ its graded 
component of weight $- \ga$ for $\ga \in \QQ_+$.
We have a canonical isomorphism
$\FF^-(\g) / \II_q^-(\g) \cong \UU^-_q(\g)$. 
By abuse of notation we will denote by the same letter 
the canonical lifting of a linear automorphism $\Phi$ of $\UU_q^-(\g)$
to an automorphism of the free algebra $\FF^-(\g)$ 
which preserves $\II_q^-(\g)$.
\medskip 
\\
\noindent
{\em{Proof of \leref{lin-le}}}. (i)
Let $\al_0 \in \chi(\Phi, \al) \cap \chi(\Phi, \al')$. 
The component of $\Phi(R_{\al, \al'})$ 
in $\FF^-_q(\g)_{- 3 \al_0}$ is 
\[
(1 - q_\al)(1- q^{-1}_\al) F_{\al_0}^3 \neq 0,
\] 
yet obviously $\II^-_q(\g)_{ - 3 \al_0} = 0$, $\forall \al_0 \in \Pi$.
This proves the first fact in part (i). 
For the second fact in (i), 
assume that $\al_1, \al_2 \in \chi(\Phi, \al)$ and 
$\al_3 \in \chi(\Phi, \al')$ for three distinct 
simple roots $\al_1, \al_2, \al_3$. The component 
of $\Phi(R_{\al, \al'})$ of weight $-\al_1 - \al_2 - \al_3$ 
is 
\[
F_{\al_3} (F_{\al_1} F_{\al_2} + F_{\al_2} F_{\al_1}) 
-(q_\al+ q^{-1}_\al) (F_{\al_1} F_{\al_3} F_{\al_2} + 
F_{\al_2} F_{\al_3} F_{\al_1}) +
(F_{\al_1} F_{\al_2} + F_{\al_2} F_{\al_1}) F_{\al_3}. 
\]
This again leads to a contradiction 
since $\II^-_q(\g)_{-\al_1 - \al_2 - \al_3}$ is 
spanned by the elements 
$[E_{\al_i}, E_{\al_j}] E_{\al_k}$ and  
$E_{\al_k} [E_{\al_i}, E_{\al_j}]$ 
for all permutations $ijk$ of $123$ such that 
$\al_i$ and $\al_j$ are not connected vertices of $\Ga$, 
(i.e., $a_{\al_i \al_j} =0$).

Part (ii) easily follows by examining in a similar way 
the nonzero elements 
$\Phi(R_{\al, \al'}) \in 
\II^-_q(\g)_{ (1- a_{\al \al'}) \theta(\al) - \theta(\al')}$
for those $\al, \al' \in \Pi$ such that $a_{\al \al'}=0$ or $-1$.
\qed
\medskip
\\
\noindent
{\em{Proof of \prref{line}}}. First, we assume that $\g$ is a 
simple Lie algebra which is not of type $B_r$ for $r \geq 3$.
In this case every linear automorphism of $\UU_q^-(\g)$
satisfies the condition in \leref{lin-le} (ii). 
If $\g$ is simply laced or of types $F_4$ or $C_r$ for $r \geq 3$, 
this follows from 
the second fact in \leref{lin-le} (i) because for these 
root systems for every simple root $\al$ of $\g$ there exists 
a simple root $\al'$ such that $a_{\al \al'} = -1$. For
root systems of rank 2 this follows from 
the first fact in \leref{lin-le} (i). Thus 
for root systems different from $B_r$, $r \geq 3$
the proposition follows from \leref{lin-le} (ii). 
  
Now, assume that $\g$ is of type $B_r$ for some $r \geq 3$. 
Denote the short simple root of $\g$ by $\al_r$ 
and the long simple roots of $\g$ by $\al_1, \ldots, \al_{r-1}$
(enumerated consecutively along $\Ga$ so that $\al_{r-1}$ 
is adjacent to $\al_r$). The second fact
in \leref{lin-le} (i) implies that there exist 
an element $\theta$ of the symmetric group $S_{\Pi}$ and  
scalars $t'_{\al_j} \in \KK^*$, $j \in [1,r]$, $t''_{\al_j} \in \KK$,
$j \in [1,r-1]$ such that
\begin{align}
\label{pr1}
\Phi(F_{\al_j}) &= t'_{\al_j} F_{\theta(\al_j)}, \; \; \forall j \in [1, r-1] \; \; 
\mbox{and} 
\\
\label{pr2}
\Phi(F_{\al_r}) &= t'_{\al_r} F_{\theta(\al_r)} 
+ \sum_{j=1}^{r-1} t''_{\al_j} F_{\theta(\al_j)}.
\end{align}
By considering the nonzero elements $\Phi(R_{\al_j, \al_{j'}}) \in 
\II^-_q(\g)_{ (1 - a_{\al_j \al_{j'}}) \theta(\al_j) - \theta(\al_{j'})}$
for $j \neq j' \in [1,r-1]$ as in \leref{lin-le} (ii)
and using that $r \geq 3$ one obtains 
$\theta(\al_r) = \al_r$. Set
$\Phi_0(F_{\al_j}) = t'_{\al_j} F_{\theta(\al_j)}$
for $j \in [1,r]$. Consider the strictly dominant integral coweight 
$\la = \sum_{j=1}^r n_j \vpi_{\al_j}\spcheck \in \PP_{++}\spcheck$
for $n_1 = \ldots = n_{r-1}=2$ and $n_r=1$.
Then $\Phi(F_{\al_j}) - \Phi_0(F_{\al_j}) \in \UU_q^-(\g)^{\geq n_j +1}$,
$\forall j \in [1,r]$
for the $\Zp$-grading of $\UU_q^-(\g)$ corresponding to $\la$, 
see \S \ref{4.1}. For graded reasons $\Phi_0$ defines
a linear automorphism of $\UU_q^-(\g)$. \leref{lin-le} (ii)
implies $\theta \in \Aut(\Ga)$ and so $\theta =\id$. Hence
$\Phi_0 = \Upsilon^{-}_{(t,\id)}$, where 
$t = (t_\al)_{\al \in \Pi}$ is given by $t_\al = (t'_\al)^{-1}$.
Eqs. \eqref{pr1}--\eqref{pr2} imply that 
$(\Upsilon^{-1}_{(t,\id)})^{-1} \Phi(F_{\al_r})$
is a $\la$-unipotent automorphism
of $\UU_q^-(\g)$. It follows from 
\thref{Uq-ut} that $\Phi = \Upsilon^-_{(t, \id)}$. 
\qed
\subsection{Proof of the main theorem}
\label{4a.3} In the remaining part of this section we will use the 
$\Zp$-grading of $\UU_q^-(\g)$ associated to 
\[
\la = \rho\spcheck = \sum_{\al \in \Pi} \vpi_\al\spcheck \in \PP_{++}\spcheck,
\]
recall \S \ref{4.1}.

The following result is due to Launois \cite[Proposition 2.3]{La3}.
We give a slightly different proof which extends to the twisted case 
under weaker assumptions on the twisting cocycle.
The lemma can be also proved using the method of 
Alev, Andruskiewitsch, and Dumas \cite[Proposition 1.2]{AAD}. 

\ble{term0} For all simple Lie algebras $\g$ of rank $r >1$, 
base fields $\KK$, $q \in \KK^*$ 
not a root of unity, and automorphisms $\Phi$ of $\UU_q(\g)$, 
we have 
\[  
\Phi(F_\al) \in \UU_q^-(\g)^{\geq 1}.
\]
\ele
\begin{proof} Launois and Lenagan proved the following fact in 
\cite[Proposition 3.2]{LaL1} (see also \cite[Corollary 2.4]{La3}):

{\em{Let $R = \oplus_{m \in \Zp} R^m$ be a connected $\Zp$-graded $\KK$-algebra 
generated in degree one by $x_i \in A^1$, $i=1,\ldots, n$ 
such that for each $i$ there exists $y_i \in R$ with
$x_i y_i = q_i y_i x_i$ for some $q_i \in \KK^*$, $q_i \neq 1$. 
Then for each automorphism 
$\Phi$ of $R$ we have $\Phi (x_i) \in R^{\geq 1}$.}}

The algebra $\UU_q^-(\g)$ satisfies the above property which 
implies the validity of the proposition. To see this, 
for each $\al \in \Pi$ choose $\al' \in \Pi$ 
such that $a_{\al \al'} \neq 0$ (recall \eqref{a}) 
and define
\[
x_{\al \al'} = \sum_{j=0}^{- a_{\al \al'}}
(-q_\al)^j \begin{bmatrix} 
-a_{\al \al'} \\ j
\end{bmatrix}_{q_\al}
      (F_\al)^j F_{\al'} (F_\al)^{-a_{\al \al'}-j}. 
\]
It follows from the quantum Serre relations \eqref{S1} that 
\[
x_{\al \al'} F_\al = q_\al^{-1} F_\al x_{\al \al'}
\]
and we have $q_\al \neq 1$.
\end{proof}
\noindent
{\em{Proof of \thref{ADc}}}. Let $\Phi \in \Aut(\UU_q^-(\g))$. 
\leref{term0} implies $\Phi(F_\al) \in \UU_q^-(\g)^{\geq 1}$, 
$\forall \al \in \Pi$. For $\al \in \Pi$, let 
$\Phi_0(F_\al) \in \UU_q^-(\g)^1 = 
\Span \{ F_{\al'} \mid \al' \in \Pi' \}$ be the unique element 
such that 
\[
\Phi(F_\al) - \Phi_0(F_\al) \in \UU_q^-(\g)^{\geq 2}.
\] 
For graded reasons, $\Phi_0$ extends to a linear
automorphism of $\UU_q^-(\g)$. Applying \prref{line} we obtain 
that there exist $\theta \in \Aut (\Ga)$ and $t \in \Tset^r$ 
such that $\Phi_0 = \Upsilon_{(t,\theta)}^-$. 
Therefore $(\Upsilon^{-}_{(t, \theta)})^{-1} \Phi$
is a $\rho\spcheck$-unipotent automorphism of $\UU_q^-(\g)$.
It follows from \thref{Uq-ut} that
$\Phi = \Upsilon_{(t, \theta)}^{-}$. This completes 
the proof of the minus case of the theorem. 
The plus case follows by applying
the isomorphism \eqref{om}.
\qed
\sectionnew{The multiparameter case}
\label{multi}
\subsection{Statement of main result}
\label{5.1}
In this section we extend \thref{ADc} to a classification 
of the automorphism groups of the 2-cocycle twists of all algebras 
$\UU^\pm_q(\g)$. This result is stated in \thref{class},
which is proved in \S \ref{5.3}. The main step is a rigidity 
result for the unipotent automorphisms of the twisted 
algebras proved in \thref{class2}. Finally, 
\S \ref{5.4} contains a classification of the isomorphisms between all 
algebras obtained by 2-cocycle twists from the algebras $\UU^\pm_{q, \pb}(\g)$
for simple Lie algebras $\g$. This is done in \thref{class}.
In particular, it is shown that 
\[
\UU_q^\pm(\g_1) \cong \UU_q^\pm(\g_2) \; \; \Leftrightarrow \; \; 
\g_1 \cong \g_2
\]
for all base fields $\KK$, non-root of unity $q$ and 
simple Lie algebras $\g_1$, $\g_2$.

Let $R$ be a $\KK$-algebra graded by an abelian group 
$C$, $R = \oplus_{c \in C} R_c$. For a 2-cocycle $\pb \in Z^2(C, \KK^*)$, 
define \cite{AST} a new algebra structure  
on the $\KK$-vector space $R$ by twisting the product in $R$ as follows:
\[
u_1 * u_2 = \pb(c_1, c_2) u_1 u_2, \; \; 
c_i \in C, u_i \in R_{c_i}, i=1,2.
\]
The twisted algebra, to be denoted by $R_\pb$, is canonically $C$-graded. 
Artin, Schelter, and Tate \cite{AST} proved that up to a graded isomorphism $R_\pb$ 
only depends on the cohomology class of $\pb$. They also proved
that, if $C$ is a free abelian group, then
\begin{equation}
\label{rb}
\rb \colon C \times C \to \KK^* \; \; 
\mbox{given by} \; \; 
\rb(c_1, c_2) = \pb(c_1, c_2) \pb(c_2, c_1)^{-1}, \; 
c_1, c_2 \in C 
\end{equation}
is a multiplicatively skew-symmetric 
group bicharacter and 
the cohomology classes $H^2(C, \KK^*)$ 
are classified by multiplicatively skew-symmetric square matrices 
of size equal to the rank of $C$ (obtained by restricting
$\rb$ to a minimal set of generators of $C$).

Given $\pb \in Z^2(\QQ, \KK^*)$, denote by $\UU^\pm_{q, \pb}(\g)$
the associated 2-cocycle twist of $\UU^\pm_q(\g)$ for the 
$\QQ$-grading from \S \ref{2.1}. The isomorphism \eqref{om}
defines an isomorphism of the twisted algebras
\begin{equation}
\label{new-om}
\om \colon \UU_{q, \pb}^\pm (\g) \to \UU_{q, \pb}^\mp(\g)
\end{equation}
because of the above mentioned property of $\rb$. The algebra 
$\UU^-_{q, \pb}(\g)$ can be described as the $\KK$-algebra
with generators $\{ F_\al \mid \al \in \Pi \}$ and 
relations
\begin{equation}
\sum_{j=0}^{1-a_{\al \al'}} (- \rb(\al', \al))^j
\begin{bmatrix} 
1-a_{\al \al'} \\ j
\end{bmatrix}_{q_\al}
      (F_\al)^j F_{\al'} (F_\al)^{1-a_{\al \al'}-j} = 0, 
\; \; \forall \al \neq \al' \in \Pi,
\label{S2}
\end{equation}
recall \eqref{S1}.

For every 2-cocycle $\pb \in Z^2(\QQ, \KK^*)$ denote by $G_\pb$ 
the subgroup of $\KK^*$ generated by the set
\[
q^2 \cup \{ q^{ - \lcor \al , \al' \rcor } \rb (\al, \al') \mid 
\al \neq \al' \in \Pi \}.
\]
If one chooses a linear ordering $<$ on $\Pi$,
then the group $G_\pb$ is also generated by the above elements 
for the pairs with $\al < \al'$.
\bde{t-f-cocyc} A 2-cocycle $\pb \in Z^2(\QQ, \KK^*)$ 
will be called torsion-free if the subgroup 
$G_\pb$ of $\KK^*$ is torsion-free.
\ede

Note that, if $\pb$ is torsion-free, then $q \in \KK^*$ is not a root of unity. 
On the other hand if $q \in \KK^*$ is not a root of unity then the 
trivial cocycle $\pb$ is torsion-free.

Denote
\[
\Aut(\Ga, \pb) = \{ \theta \in \Aut(\Ga) \mid 
\rb(\theta(\al), \theta(\al') ) = 
\rb(\al, \al'), \; \; 
\forall \al, \al' \in \Pi \}.
\]
We have an embedding 
$\Upsilon^\pm \colon \Tset^r \rtimes \Aut (\Ga, \pb) \hra \Aut( \UU_{q, \pb}^\pm(\g))$,
where for $(t, \theta) \in \Tset^r \rtimes \Aut (\Ga, \pb)$ the automorphism
$\Upsilon_{(t, \theta)}^\pm \in \Aut( \UU_{q, \pb}^\pm(\g))$ is 
given by \eqref{Upsilon1}--\eqref{Upsilon2}.

\bth{class} For every simple Lie algebra $\g$ of rank $r >1$,
base field $\KK$, $q \in \KK^*$, 
and a torsion-free 2-cocycle $\pb \in Z^2(\QQ, \KK^*)$ 
satisfying
\begin{equation}
\label{cond2}
q_\al \rb(\al,\al'), q_\al^{-1} \rb(\al, \al') \neq 1, \; 
\forall \al, \al' \in \Pi \; \; \mbox{such that} \; \;  
a_{\al \al'} = - 1,
\end{equation}
the map
\[
\Upsilon^\pm \colon \Tset^r \rtimes \Aut (\Ga, \pb) \to \Aut( \UU_{q, \pb}^\pm(\g))
\] 
is a group isomorphism.
\eth

The special case of $\g = \so_5$ of this theorem was obtained by Tang \cite{T}.
Because of the isomorphism \eqref{new-om}, it is sufficient to prove the 
theorem in the minus case.

We finish this subsection with a result which explains the origin of the 
torsion-free condition from \deref{t-f-cocyc}. Let $\ib= (\al_1, \ldots, \al_N)$ be 
a reduced word for the longest element $w_0$ of $W$. All automorphisms 
and skew derivations in the iterated Ore extension presentation 
\eqref{Umin-iter} are graded. Thus for all $\pb \in Z^2(\QQ, \KK^*)$,
we have the iterated Ore extension presentation 
\begin{equation}
\label{Umin-iter-twist}
\UU_{q, \pb}^-(\g) = \KK [F_{\be_1}] [F_{\be_2}; \sig_2, \delta_2] \ldots 
[F_{\be_N}; \sig_N, \delta_N ],
\end{equation}
where $\sig_l$ and $\delta_l$ are still given by \eqref{sig} and 
\eqref{delta} but this time $t_l \in \Tset^r$ are such that 
$t_l^{\be_k} = q^{ \lcor \be_l , \be_k \rcor } \rb(\be_l, \be_k)^{-1}$,
$\forall k \in [1,l]$. If $q \in \KK^*$ is not a root of unity,
this is a CGL extension for the following 
choice of the elements $q_{lk}$ and $q_l \in \KK^*$ (recall \deref{CGL}):
\begin{equation}
\label{q-elem}
q_{lk} = q^{ - \lcor \be_l , \be_k \rcor } \rb(\be_l, \be_k), \; 
1 \leq k < l \leq N \; \; 
\mbox{and} \; \; 
q_l = q_{\al_l}^{-2}, \; l \in [1,N].
\end{equation}

\bpr{t-f-cocy} Let $\g$ be a simple Lie algebra of rank $r >1$.
For all 2-cocycles $\pb \in Z^2(\QQ, \KK^*)$ and 
reduced words $\ib$ for the longest element $w_0$ of $W$, the group 
$G_\pb$ is precisely the subgroup of $\KK^*$ generated 
by the elements $q_{lk} \in \KK^*$, $1 \leq k < l \leq N$ 
given by \eqref{q-elem}.

In particular, $\pb \in Z^2(\QQ, \KK^*)$ is a torsion-free cocycle
if and only if the iterated Ore extension presentation 
\eqref{Umin-iter-twist} of $\UU_{q, \pb}^-(\g)$ associated to 
one reduced word $\ib$ for $w_0$ (and thus to every
reduced word $\ib$ for $w_0$) is  
is a torsion-free CGL extension.
\epr
\begin{proof} 
Denote by $G_\ib$ the subgroup of $\KK^*$ generated by the elements 
$q_{lk}$, $1 \leq k < l \leq N$ in Eq. \eqref{q-elem}. 
First we show that $G_\pb \subseteq G_\ib$. 
Assume that $\al, \al' \in \Pi$ and $\al$ comes before $\al'$ in the 
ordering 
\begin{equation}
\label{list}
\be_1, \ldots, \be_N
\end{equation}
from Eq. \eqref{rootv}
of the positive roots of $\g$. If $\al$ and $\al'$ are 
not connected with an edge in $\Ga$, 
then $\lcor \al, \al' \rcor = 0$ and 
$q^{ - \lcor \al, \al' \rcor } \rb(\al, \al') = 
(q^{ - \lcor \al', \al \rcor } \rb(\al', \al) )^{-1}
\in G_\ib$. If they are connected by an edge, then 
the root $\al + \al'$ of $\g$ is listed between 
$\al$ and $\al'$ in \eqref{list} since the 
ordering \eqref{list} is convex and
\begin{align*}
q^{ - \lcor \al, \al \rcor } &= 
( q^{ -\lcor \al + \al', \al \rcor} \rb (\al + \al', \al) ) 
( q^{ - \lcor \al', \al \rcor} \rb (\al', \al) )^{-1} \in G_\ib,
\\
q^{ - \lcor \al', \al' \rcor } &= 
( q^{ -\lcor \al', \al + \al' \rcor} \rb (\al', \al + \al') ) 
( q^{ - \lcor \al', \al \rcor} \rb (\al', \al) )^{-1} \in G_\ib.
\end{align*}
Thus $q^2 \in G_\ib$ and $q^{ - \lcor \al, \al' \rcor } \rb(\al, \al') = 
q^{- 2 \lcor \al, \al' \rcor}
(q^{ - \lcor \al', \al \rcor } \rb(\al', \al) )^{-1} \in G_\ib$.
Hence $G_\pb \subseteq G_\ib$. For the opposite inclusion, 
we fix a linear ordering $<$ on $\Pi$. If
$\be_l = \sum_{\al \in \Pi} m_\al \al$ and 
$\be_k = \sum_{\al \in \Pi} n_\al \al$ for $k < l$, then
\[
q^{- \lcor \be_l, \be_k \rcor} \rb (\be_l, \be_k) = 
\prod_{\al \in \Pi} q^{ - m_\al n_\al \lcor \al, \al \rcor}  
\prod_{ \al < \al' \in \Pi}
q^{ - 2 m_{\al'} n_\al \lcor \al, \al' \rcor } 
\left( q^{-\lcor \al, \al' \rcor } \rb(\al, \al') \right)^{ 
m_\al n_{\al'} - m_{\al'} n_\al} \in G_\pb.
\] 
This completes the proof of the proposition.
\end{proof}

The above argument shows that the group $G_\pb$ can be also characterized
as the subgroup of $\KK^*$ generated by all elements of the 
form $q^{- \lcor \be, \be' \rcor } \rb(\be, \be')$, where 
$\be$ and $\be'$ run over all positive roots of $\g$.
\subsection{Unipotent automorphisms of $\UU_{q, \pb}^-(\g)$} 
\label{5.2}
Each strictly dominant integral coweight $\la= \sum_{\al \in \Pi} m_\al \vpi_\al\spcheck$
gives rise to a specialization of the $(-\QQ_+)$-grading
of $\UU_{q, \pb}^-(\g)$ to a connected $\Zp$-grading
by setting $\deg F_\al := m_\al = \lcor \la, \al \rcor$
for $\al \in \Pi$. We will 
denote the corresponding graded components by $\UU_{q, \pb}^-(\g)^m$, 
$m \in \Zp$. Analogously to the untwisted case we will 
call an automorphism $\Phi$ of $\UU_{q, \pb}^-(\g)$
$\la$-unipotent if
\[
\Phi(F_\al) - F_\al \in \UU_{q, \pb}^-(\g)^{\geq \lcor \la, \al \rcor + 1},
\; \; 
\forall \al \in \Pi.
\]

\bth{class2} Let $\g$ be a Lie algebra of rank $r >1$,
$\KK$ and arbitrary base field, $q \in \KK^*$, and $\pb \in Z^2(\QQ, \KK^*)$
a torsion-free 2-cocycle. Every $\la$-uniponet 
automorphism $\Phi$ of $\UU_{q, \pb}^-(\g)$ 
for a strictly dominant integral coweight $\la$
is equal to the identity automorphism.
\eth
\begin{proof}
Let $\ib$ be a reduced word for $w_0$. It follows from 
\prref{t-f-cocy} that \eqref{Umin-iter-twist} 
is a torsion--free CGL extension
presentation of $\UU_{q, \pb}^-(\g)$.
Denote by $(\ol{F}_{\ib, 1}, \ldots, \ol{F}_{\ib, N})$ 
the final $N$-tuple from the Cauchon deleting derivation procedure 
applied to it. Let $\qb$ denote the multiplicatively skew-symmetric 
$N \times N$ matrix whose entries $q_{lk}$, $1 \leq k < l \leq N$ 
satisfy eq. \eqref{q-elem}. Then we have an isomorphism 
of quantum tori 
\[
\TT_\qb \cong \TT(\ib, \pb) := \KK \lcor \ol{F}_{\ib, 1}^{ \, \pm 1}, 
\ldots \ol{F}_{\ib, N}^{ \, \pm 1} \rcor
\subseteq \Fract(\UU_{q, \pb}^-(\g)) \; \; 
\mbox{given by} \; \; 
X_l \mt \ol{F}_{\ib, l}, l \in [1,N]
\]  
recall \eqref{q-tor}. Since the CGL extension presentation 
\eqref{Umin-iter-twist} is torsion-free, 
the quantum torus $\TT_\qb$ is saturated, see \S \ref{2.2}.
Eq. \eqref{main1b-eq} 
is a graded equality in $\UU_q^-(\g)$ and thus it holds in 
$\UU_{q, \pb}^-(\g)$ after an appropriate rescaling. Thus 
\[
\De_{\ib, 1}, \ldots, \De_{\ib, N}
\]
is a generating set of the quantum torus $\TT(\ib, \pb)$. 
Recall from \S \ref{3.1} that the property of a quantum torus 
being saturated does not depend on the choice of its 
generators.
We use the degree
vector $\db$ from \eqref{deg-vect} to define a $\Zset$-grading on $\TT(\ib, \pb)$ and 
to form a completion as in \S \ref{3.1}. This completion will be denoted by 
$\wh{\TT}(\ib, \pb, \db)$. Analogously to \S \ref{4.1}, 
to every $\la$-unipotent automorphism $\Phi$ of $\UU_{q, \pb}^-(\g)$
we associate a bi-finite unipotent automorphism 
of the completed saturated quantum torus $\wh{\TT}(\ib,\pb, \db)$.

Analogs of 
Propositions \ref{pste1} and \ref{p2re} hold under very mild modifications.
Denote by $\NN_{q, \pb}^-(\g)$ the subalgebra of $\UU_{q, \pb}^-(\g)$ 
generated (and spanned) by $b^\la_1$, $\la \in \PP_+$. It is a 
quantum affine space algebra with generators $b^{\vpi_\al}_1$, $\al \in \Pi$ 
and relations $b^{\vpi_\al}_1 b^{\vpi_{\al'}}_1 = 
\rb( (1-w_0) \vpi_\al, (1-w_0)\vpi_{\al'}) b^{\vpi_{\al'}}_1 b^{\vpi_\al}_1$,
$\forall \al, \al' \in \Pi$, recall \eqref{norm}. Consider 
the localization $\NN_{q, \pb}^-(\g)^\sharp := \NN_{q,\pb}^-(\g) [\Om(\g)^{-1}]$ 
where $\Om(\g) = \{ b_1^\la \mid \la \in \PP_+ \}$
and the elements $b^\mu_1$, $\mu \in \PP$ in it given by \eqref{b-mu-elem}.
For all $\mu \in \PP$, we have 
\[
b_1^\mu u
= 
q^{\lcor (1+ w_0)\mu, \ga \rcor }
\rb( (1-w_0) \mu, \ga)
u b_1^\mu, \; \; 
\forall \, 
u \in \UU_q^-(\g)_{-\ga},
\ga \in \QQ_+. 
\]
This property and the first part of \thref{1a} (ii) imply
\[
Z ( \UU_{q, \pb}^-(\g) [\Om(\g)^{-1}]) 
= \CC_{q, \pb}^-(\g)^\sharp:= 
\{ b^\mu_1 \mid \mu \in \PP, \; 
q^{\lcor (1+ w_0)\mu, \ga \rcor }
\rb( (1-w_0) \mu, \ga)=1, \; \forall 
\ga \in \QQ \}.
\]
The argument of the proof of Eq. \eqref{ZZZ} gives
\[
Z( \TT(\ib, \pb)) = \CC_{q, \pb}^-(\g)^\sharp.
\]
For graded reasons it 
follows from \thref{1a} (iii) that $\UU_{q, \pb}^-(\g)$ 
is a free left and right $\NN_{q, \pb}^-(\g)$-module 
with basis \eqref{basis}. Analogously to the proofs of 
Lemmas \ref{l43a} and \ref{ldir-sum} and \prref{2re}
we obtain
\[
Z( \UU_{q, \pb}^-(\g) ) 
= \CC_{q, \pb}^-(\g):= \Span
\{ b^\mu_1 \mid \mu \in \PP_+, \; 
q^{\lcor (1+ w_0)\mu, \ga \rcor }
\rb( (1-w_0) \mu, \ga)=1, \; \forall 
\ga \in \QQ \}
\]
and then show that for every $\la$-unipotent 
automorphism $\Phi \in \Aut (\UU_{q, \pb}^-(\g))$ 
there exist elements
$z_\al \in \CC_{q, \pb}^-(\g)^{\geq 1}$ for $\al \in \Pi$ such that
\begin{equation}
\label{z-elem}
\Phi(F_\al) = (1+ z_\al) F_\al, \; \; 
\forall \al \in \Pi.
\end{equation}
Recall \thref{1a} (i).
The argument of the proof of \cite[Theorem 4.1]{GLen-qd} 
of Goodearl and Lenagan
shows that $\UU_{q, \pb}^-(\g) b^{\vpi_\al}_1$ 
are height one prime ideals of $\UU_{q, \pb}^-(\g)$ for all $\al \in \Pi$. 
Using this, analogously to \S \ref{4.4} one shows that 
$z_\al =0$, $\forall \al \in \Pi$. This completes
the proof of the theorem.
\end{proof}
\subsection{Proof of \thref{class}}
\label{5.3}
To each semisimple Lie algebra $\g$ one can attach
a $\KK$-algebra $\UU^-_q(\g)$ analogously to 
\S \ref{2.1} using the normalized $W$-invariant bilinear 
form on $\Rset \Pi$ such that $\lcor \al, \al \rcor = 2$ 
for all short simple roots $\al$ of $\g$. We will need those
algebras for induction purposes.
For a subset $\Pi' \subset \Pi$ and $\pb \in Z^2(\QQ, \KK^*)$,
denote by $\UU_{q, \pb}^-(\g_{\Pi'})$
the subalgebra of $\UU_{q, \pb}^-(\g)$ generated by $F_\al$ for 
$\al \in \Pi'$. An algebra automorphism
$\Phi \in \Aut ( \UU_{q, \pb}^-(\g_{\Pi'}) )$
will be called linear if 
$\Phi (F_\al) \subseteq \Span \{ F_{\al'} \mid \al' \in \Pi' \}$,
$\forall \al \in \Pi'$. We will use the notation from Eq. \eqref{chi} for those.

Given $\pb \in Z^2(\QQ, \KK^*)$, for $c \in \KK^*$ denote 
\[
\Pi^c = \{ \al \in \Pi \mid \exists \al' \in \Pi \; \; 
\mbox{such that} \; \; a_{\al \al'} = 0 \; \; 
\mbox{and} \; \; \rb(\al', \al) = c \},
\]
recall \eqref{a}. Denote 
\[
\Pi^* = \cup_{c \in \KK^*, c \neq 1} \Pi^c, \; \; 
\mbox{and} \; \; 
\Pi^\bullet = \Pi \backslash \Pi^*.
\]

The following lemma is proved analogously 
to \leref{lin-le}. 
\ble{lin2} For all semisimple Lie algebras $\g$, base fields $\KK$, 
$q \in \KK^*$, $\pb \in Z^2(\QQ, \KK^*)$, and linear automorphisms
$\Phi$ of $\UU_{q, \pb}^-(\g)$ the following hold:

(i) If $\al, \al' \in \Pi$ are such that $a_{\al \al'} = 0$ 
and $\rb(\al', \al) \neq 1$ then 
$\chi(\Phi, \al) \cap \chi(\Phi, \al') = \emptyset$.

(ii) If $\al \in \Pi^c$ for some $c \neq 1$, 
then $\chi(\al, \Phi) \subseteq \Pi^c$.  

(iii) If the condition \eqref{cond2} is satisfied, then 
for all $\al, \al' \in \Pi$ such that 
$a_{\al \al'} = -1$ we have 
$\chi(\Phi, \al) \cap \chi(\Phi, \al') = \emptyset$. 
If, in addition, for such a pair $(\al, \al')$
we have $\chi(\Phi, \al) \subseteq \Pi^\bullet$, then 
$|\chi(\Phi, \al)| = 1$.

(iv) If there exist an element $\theta$ of the symmetric group $S_{\Pi}$ and  
scalars $t'_\al \in \KK^*$ for $\al \in \Pi$ such that
\[
\Phi(F_\al) = t'_\al F_{\theta(\al)}, \; \; \forall \al \in \Pi,
\]
then $\theta \in \Aut(\Ga, \pb)$ and $\Phi = \Upsilon^-_{(t,\theta)}$,
where $t= (t_\al)_{\al \in \Pi} \in \Tset^r$ is given by
$t_\al = t_{\theta^{-1}(\al)}^{-1}$. 
\ele
Next we prove an extension of \prref{line} to the twisted case.
\bpr{lin-twist} Let $\g$ be a semisimple Lie algebra of rank $r$, 
$\KK$ an arbitrary base field, $q \in \KK^*$, and 
$\pb \in Z^2(\QQ, \KK^*)$ a torsion-free 2-cocycle satisfying \eqref{cond2}.
Then every linear automorphism $\Phi$ of $\UU_{q, \pb}^-(\g)$ 
is of the form $\Upsilon_{(t, \theta)}^-$ for some 
$\theta \in \Aut(\Ga, \pb)$ and $t \in \Tset^r$. 
\epr
\begin{proof} First, we claim that for every $\al \in \Pi^*$ 
there exists a subset $\Pi' \subseteq \Pi^*$ containing $\al$ 
such that $\Phi$ restricts to a linear automorphism of 
$\UU_{q, \pb}^-(\g_{\Pi'})$ and 
\begin{equation}
\label{co11}
\rb(\al_1, \al_2) = 1, \; \; 
\forall \al_1, \al_2 \in \Pi' \; \; 
\mbox{such that} \; \; a_{\al_1 \al_2} = 0.
\end{equation}
Choose $c \in \KK^*$, $c \neq 1$ such that $\al \in \Pi^c$.
By \leref{lin2} (ii), $\Phi$ restricts to a linear 
automorphism of $\UU_{q, \pb}^-(\g_{\Pi^c})$. 
If $\Pi'=\Pi^c$ satisfies \eqref{co11}, this proves the 
claim. Otherwise we continue recursively
by using $\g_{\Pi^c}$ in place of $\g$. Analogously 
to the proof of \prref{line}, the claim 
and \leref{lin2} (iii) imply that 
there exist $\theta^* \in S_{\Pi^*}$ 
and $t^* = (t'_\al)_{\al \in \Pi^*} \in \Tset^{|\Pi^*|}$ such that 
\begin{equation}
\label{ll5}
\Phi(F_\al) = t'_\al F_{\theta^*(\al)}, \; \; 
\forall \al \in \Pi^*.
\end{equation}
Let
\[
\Phi(F_\al) = \sum_{\al' \in \Pi^\bullet} c_{\al \al'} F_{\al'} + 
\sum_{\al'' \in \Pi^*} c_{\al \al''} F_{\al''}, \; \;  
\al \in \Pi^\bullet.
\]
for some $c_{\al \al'}, c_{\al \al''} \in \KK$. 
It follows from \eqref{ll5} and the form of the quantum Serre relation \eqref{S2} that 
\[
\Phi^\bullet(F_\al) := \sum_{\al' \in \Pi^\bullet} c_{\al \al'} F_{\al'}, 
\; \; \al \in \Pi^\bullet
\]
extends to a linear automorphism of $\UU_{q, \pb}^-(\g_{\Pi^\bullet})$.
Analogously to the proof of \prref{line}, using 
\leref{lin2} (iii) one obtains that there exist
$\theta^\bullet \in S_{\Pi^\bullet}$ and $t^\bullet = (t''_\al) \in 
\Tset^{|\Pi^\bullet|}$ such that 
\[
\Phi^\bullet(F_\al) = t''_\al F_{\theta^\bullet(\al)}, \; \; 
\forall \al \in \Pi^\bullet.
\]
Consider the $\Zp$-grading of $\UU_{q, \pb}^-(\g)$ associated 
to the strictly dominant integral coweight 
$\la = \sum_{\al \in \Pi} n_\al \vpi_{\al}\spcheck$, 
where $n_\al = 1$ if $\al \in \Pi^\bullet$ and $n_\al = 2$
if $\al \in \Pi^*$. For graded reasons it follows that 
\[
\Phi_0(F_\al) := 
\begin{cases}
t'_\al F_{\theta^\bullet(\al)}, & \mbox{if} \; \; \al \in \Pi^\bullet
\\
t''_\al F_{\theta^*(\al)}, & \mbox{if} \; \; \al \in \Pi^*
\end{cases}
\]
extends to a linear automorphism of $\UU_{q, \pb}^-(\g)$ 
such that 
\begin{equation}
\label{ll6}
\Phi(F_\al) - \Phi_0(F_\al) \in \UU_{q, \pb}^-(\g)^{\geq n_\al + 1}, \; 
\forall \al \in \Pi. 
\end{equation}
\leref{lin2} (iv) implies that $\Phi_0 = \Upsilon^-_{(t,\theta)}$ for 
some $\theta \in \Aut(\Ga, \pb)$ and $t \in \Tset^r$. 
By Eq. \eqref{ll6}, $(\Upsilon_{(t, \theta)}^-)^{-1} \Phi$ 
is a $\la$-unipotent automorphism of $\UU_{q, \pb}^-(\g)$ and 
by \thref{class2} $\Phi = \Upsilon_{(t, \theta)}^-$. 
\end{proof}
The last step before the proof of \thref{class} is an 
extension of \leref{term0} to the twisted case.

\ble{term1} Let $\g$ be a simple Lie algebra, $\KK$ an
arbitrary base field, $q \in \KK^*$ not a root of unity,
and $\pb \in Z^2(\QQ, \KK^*)$. For all automorphisms
$\Phi$ of $\UU_{q,\pb}^-(\g)$ we have 
\[
\Phi(F_\al) \in \UU_{q, \pb}^-(\g)^{\geq 1}
\]
with respect to the $\Zp$-grading of $\UU_{q, \pb}^-(\g)$ 
associated to $\la = \rho\spcheck = \sum_{\al \in \Pi} \vpi_\al\spcheck$.
\ele
\begin{proof} For $\al \neq \al'$ denote
\[
x_{\al \al'}^\pm = \sum_{j=0}^{- a_{\al \al'}}
(- \rb(\al',\al) q_\al^{\pm 1})^j \begin{bmatrix} 
-a_{\al \al'} \\ j
\end{bmatrix}_{q_\al^{\pm 1} }
      (F_\al)^j F_{\al'} (F_\al)^{-a_{\al \al'}-j}. 
\]
It follows from \eqref{S2} that 
\[
x_{\al \al'}^\pm F_\al = \rb (\al', \al) q_\al^{\mp a_{\al \al'}} F_\al x_{\al \al'}.
\]
If $a_{\al \al'} \neq 0$, then either $\rb (\al', \al) q_\al^{a_{\al \al'}} \neq 1$ or 
$\rb (\al', \al) q_\al^{-a_{\al \al'}} = 1$ because $q$ is not a root of unity.
This establishes the lemma analogously to the proof of \leref{term0}.
\end{proof}
\noindent
{\em{Proof of \thref{class}}}. Because of the isomorphism from Eq. \eqref{new-om}
it is sufficient to prove the minus case of the theorem. \leref{term1} implies that 
every automorphism $\Phi \in \Aut(\UU_{q, \pb}^-(\g))$ 
satisfies $\Phi(F_\al) \in \UU_{q, \pb}^-(\g)^{\geq 1}$ 
with respect to $\Zp$-grading from \S \ref{5.2} 
corresponding to the strictly dominant integral coweight 
$\la = \rho\spcheck = \sum_{\al \in \Pi} \vpi_\al\spcheck$. 
For each $\al \in \Pi$ 
there exists a unique element $\Phi_0(F_\al) \in \UU_{q, \pb}^-(\g)^1$ 
such that 
\begin{equation}
\label{PhiPhi0}
\Phi(F_\al) - \Phi_0(F_\al) \in \UU_{q, \pb}^-(\g)^{\geq 2}.
\end{equation}
For graded reasons $\Phi_0$ extends to a linear automorphism of 
$\UU_{q, \pb}^-(\g)$. Hence, by \prref{lin-twist}, $\Phi_0 = \Upsilon_{(t, \theta)}^-$
for some $\theta \in \Aut(\Ga, \pb)$ and $t \in \Tset^r$. Moreover, it follows from
\eqref{PhiPhi0} that $\Phi_0^{-1} \Phi$ is a $\rho\spcheck$-unipotent 
automorphism of $\UU_{q, \pb}^-(\g)$. \thref{class2} implies that 
$\Phi = \Phi_0$ and thus $\Phi =\Upsilon_{(t, \theta)}^-$. 
\qed
\bre{semi}
It is easy to see that the proof of \thref{class} easily extends to the more general statement 
for all semisimple Lie algebras $\g$ not having a simple direct summand of rank 1. 
For simplicity of the exposition we only considered the simple case. The details of the 
semisimple case are left to the reader. 
\ere
\subsection{The isomorphism problem for the algebras $\UU_{q, \pb}^\pm(\g)$}
\label{5.4}
The following theorem applies the results on automorphism groups from this and the previous 
sections to the 
isomorphism problem for the class of algebras of the form $\UU^\pm_{q, \pb}(\g)$
for simple Lie algebras $\g$.
The idea to apply results on automorphism groups to this isomorphism problem 
was suggested to us by Len Scott.

\bth{isom} Let $\g_i$, $i=1,2$ be two simple Lie algebras 
with root lattices $\QQ_i$, Dynkin diagrams $\Ga_i$
and set of simple roots $\Pi_i$ (considered as the 
set of vertices of $\Ga_i$). Let $\KK$ be 
an arbitrary base field, $q \in \KK^*$, and $\pb_i \in Z^2(\QQ_i, \KK^*)$
be two 2-cocycles satisfying \eqref{cond2} and such that 
$G_{\pb_1} G_{\pb_2}$ is a torsion-free subgroup of $\KK^*$. 
Then 
\begin{equation}
\label{is1}
\UU_{q, \pb_1}^\pm(\g_1) \cong \UU_{q, \pb_2}^\pm(\g_2)
\end{equation}
if and only if 
\begin{align}
\label{is2}
&\g_1 \cong \g_2 \; \; \mbox{and there exists a 
graph isomorphism} \; \; \theta \colon \Ga_1 \to \Ga_2
\\
&\mbox{such that} \; \; 
\rb_2(\theta(\al), \theta(\al')) = \rb_1(\al, \al'), \; \; 
\forall \al, \al' \in \Pi_1.
\nn
\end{align}
Here, as in \eqref{rb}, we set $\rb_i(\al, \al'): = \pb(\al, \al') \pb(\al', \al)^{-1}$
for $\al , \al' \in \Pi_i$.
\eth
\bco{isom2} For all simple Lie algebras $\g_1$ and $\g_2$, base fields
$\KK$ and a non-root of unity $q \in \KK^*$,
\[
\UU_q^\pm(\g_1) \cong \UU_q^\pm(\g_2) 
\quad \mbox{if and only if} \quad
\g_1 \cong \g_2.
\]
\eco
\noindent
{\em{Proof of \thref{isom}.}} It is obvious that \eqref{is2} implies \eqref{is1}. 
The Gelfand--Kirillov dimensions of $\UU_{q, \pb}(\g_i)$ 
are equal to the number of positive roots of $\g_i$ 
because of the iterated Ore extension presentation \eqref{Umin-iter-twist}.
Therefore \eqref{is1} $\Rightarrow$ \eqref{is2} in the case 
when the rank of one of the algebras $\g_i$ is equal to $1$.
Thus all we need to show is that 
\[
\mbox{\eqref{is1}} \Rightarrow \mbox{\eqref{is2}} 
\; \; \mbox{in the case when $\rank(\g_1), \rank(\g_2)>1$}.
\]
Because of the isomorphism \eqref{new-om} we can restrict 
to the minus case. 

Set $\g = \g_1 \oplus \g_2$. 
Its root lattice, Dynkin diagram and set of simple roots 
are given by $\QQ = \QQ_1 \oplus \QQ_2$, $\Ga = \Ga_1 \sqcup \Ga_2$
and $\Pi = \Pi_1 \sqcup \Pi_2$. There is a unique cocycle 
$\pb \in Z^2(\QQ, \KK^*)$ such that
$\pb(\al, \al') := \pb_i(\al, \al')$ if $\al, \al' \in \QQ_i$
for some $i=1,2$ 
and $\pb(\al, \al') =1$ if $\al \in \QQ_1, \al' \in \QQ_2$
or $\al \in \QQ_2, \al' \in \QQ_1$. 

The conditions on $\pb_1$ and $\pb_2$ 
are equivalent to saying that $\pb$ is a torsion-free 2-cocycle 
that satisfies \eqref{cond2}. We also have a canonical isomorphism 
\[
\UU_{q, \pb}^-(\g) \cong \UU_{q, \pb_1}^-(\g_1) \otimes \UU_{q, \pb_2}^-(\g_2).
\]

Let $\Psi \colon \UU_{q, \pb_1}^-(\g_1) \to \UU_{q, \pb_2}^-(\g_2)$
be an algebra isomorphism. 
Following Kimmerle \cite{K}, consider the involutive automorphism
\begin{equation}
\label{Phi}
\Phi \in \Aut(\UU_{q, \pb}^-(\g)), \quad
\mbox{given by} \; \; 
\Phi(u_1 \otimes u_2) := \Psi^{-1}(u_2) \otimes \Psi(u_1), \; \; 
\forall u_i \in \UU_{q, \pb_i}(\g_i).
\end{equation}
It is easy to see that the statement of \leref{term1} holds for every
semisimple Lie algebra $\g$ without simple direct summands of rank 1.
It follows from \leref{term1} that $\Phi(F_\al) \in \UU_{q, \pb}^-(\g)^{\geq 1}$ 
$\forall \al \in \Pi$ with respect to $\Zp$-grading from \S \ref{5.2} 
corresponding to  
$\la = \rho\spcheck = \sum_{\al \in \Pi} \vpi_\al\spcheck$. 
As in the proof of \thref{class}, for each $\al \in \Pi$ 
there exists a unique element 
$\Phi_0(F_\al) \in \UU_{q, \pb}^-(\g)^1$ such that 
\[
\Phi(F_\al) - \Phi_0(F_\al) \in \UU_{q, \pb}^-(\g)^{\geq 2}.
\]
For graded reasons $\Phi_0$ extends to a linear automorphism of 
$\UU_{q, \pb}^-(\g)$. Applying \prref{lin-twist}, 
we obtain that $\Phi_0 = \Upsilon_{(t, \theta_0)}^-$ for some $\theta_0 \in \Aut(\Ga, \pb)$ 
and $t \in \Tset^{|\Pi|}$. The definition \eqref{Phi} of $\Phi$ 
implies that $\theta_0(\Ga_1) = \Ga_2$ and $\theta_0(\Ga_2) = \Ga_1$. 
From the definitions of $\Aut(\Ga, \pb)$ and $\pb$ we obtain
\[
\rb_2(\theta_0(\al), \theta_0(\al')) = \rb_1(\al, \al'), \; \; 
\forall \al, \al' \in \Pi_1.
\]
Thus \eqref{is2} holds for $\theta:= \theta_0|_{\Ga_1}$,
which completes the proof of the theorem.
\qed

\end{document}